\documentclass[12pt,draft,leqno]{article}
\usepackage{amssymb,eucal,latexsym}
\usepackage{amsmath,amsthm,latexsym,amscd,amsbsy,fleqn,graphicx}
\usepackage{amsmath,amssymb}

\usepackage{amsthm}
\usepackage[T2A]{fontenc}
\usepackage[cp1251]{inputenc}
\usepackage[english, bulgarian]{babel}

\usepackage{amsthm}

\newtheorem{theorem}{Theorem}[subsection] \newtheorem{lm}[theorem]{Lemma} \newtheorem{exa}[theorem]{Example} \newtheorem{cor}[theorem]{Corollary} \newtheorem{pro}[theorem]{Proposition}
\newtheorem{defi}[theorem]{Definition} \newtheorem{nota}[theorem]{Notation} \newtheorem{rem}[theorem]{Remark} \newtheorem{nist}[theorem]{}

\def\unlx{\underline{X}}
\def\Set{\mbox{{\boldmath $Set$}}} \def\CohSp{\mbox{{\boldmath $CohSp$}}} \def\DLat{\mbox{{\boldmath $DLat$}}} \def\OStone{\mbox{{\boldmath $OStone$}}} \def\T{{\cal
T}} \def\L{{\cal L}} \def\PP{{\cal P}} \def\BB{{\cal B}} \def\HH{{\cal H}} \def\ion{i=1,\ldots,n} \def\jom{j=1,\ldots,m}

\def\op{\oplus} \def\ot{\otimes} \def\we{\wedge} 

\def\p{\varphi} \def\a{\alpha} \def\b{\beta}   \def\OM{\Omega}  \def\s{\sigma} \def\t{\theta} \def\TE{\Theta}

\def\sq{\hbox{\vrule\vbox{\hrule\phantom{o}\hrule}\vrule}}  \def\R{\mbox{{\boldmath $R$}}} 
 \def\D{\mbox{{\boldmath $D$}}}   \def\SS{\mbox{{\boldmath $S$}}}

\def\Ra{\Rightarrow} \def\La{\Leftarrow} \def\lra{\longrightarrow} \def\sbe{\subseteq}   \def\stm{\setminus}
\def\ems{\emptyset}     \def\doc{\hspace{-1cm}{\it Proof.}~~}

\title{{\huge\bf On the Investigations}\\
\vspace{0.2cm}
{\huge\bf of Ivan Prodanov}\\
\vspace{0.2cm}
{\huge\bf in the Theory of Abstract Spectra}\\
\vspace{0.4cm}
{\large\bf Georgi Dimov and Dimiter Vakarelov }\\
\vspace{0.1cm}
{\footnotesize\rm Department of Mathematics and Informatics,  University of Sofia,}\\
{\footnotesize\rm 5 J. Bourchier Blvd., 1164 Sofia, Bulgaria}\\
\vspace{0.2cm}
{\Large\it Dedicated to the memory of Ivan
Prodanov}\\
}

\author{}

\date{}

\begin{document}

\maketitle

\baselineskip = 0.8\normalbaselineskip


\begin{center}
{\footnotesize Abstract}
\end{center}

\smallskip

{\footnotesize
We were invited by the Organizing Committee of the Mathematical Conference dedica\-ted to Professor Ivan Prodanov on the occasion of the 60-th anniversary of his birth and the 10-th anniversary of his death, held on May 16, 1995 at the Faculty of
Mathematics and Informatics of the Sofia University $``$St. Kliment Ohridski'', to give a talk on his investigations in the theory of abstract spectra. All of his results in this area were announced by him
in a short paper published in \cite{P1}
and, as far as we know, their proofs were never written by him in the form of a manuscript, preprint or paper. The very incomplete notes
which we have from Prodanov's talks on the Seminar on Spectra, organized by him in 1979/80 academic year at the Faculty of Mathematics of the University of Sofia,
seem to be the only trace of a small part of the proofs of some of the results from the cited above paper. Since the untimely death of Ivan Prodanov withheld him from preparing the
full version of this paper and since, in our opinion, it contains interesting and important results, we undertook the task of writing a full
version of it and thus making the results from it known to the mathematical community. So, the aim of this paper is to supply with proofs the results of Ivan Prodanov announced in the cited above paper, but we added also a small amount of new results. The main of them is Theorem \ref{th66} which was formulated and proved by us as a generalization of Prodanov's assertions Corollary \ref{th6} and Corollary \ref{th7}.

The full responsibility for the correctness of the proofs of the assertions presented below in this work is taken by us;  just for this reason our names appear as
authors of the present paper.
The talks of the participants of the  conference  had to be published at a special volume of the Annuaire de l’Universite de Sofia $``$St. Kliment Ohridski”,
but this never happened. That is why we have decided to publish our work separately. Since our files were lost and we had to write them once more,  the paper appears only now.
}

\smallskip


\noindent{\footnotesize {\em  2010 MSC:} Primary 06D05; Secondary: 54H10, 54B30, 54E55, 06D50, 18A40, 18B30.

\smallskip

\noindent{\em Keywords:} Abstract spectra, bitopological spaces, compact spaces,  coherent spaces, coherent maps, distributive lattices, (pre)separative algebras, convex space, separation theorems.}

\footnotetext[1]{{\footnotesize {\em E-mail addresses:}
gdimov@fmi.uni-sofia.bg, dvak@fmi.uni-sofia.bg}}

\baselineskip = 1.5\normalbaselineskip


\section{Introduction}

\hspace{1cm}This paper contains an extended version of the invited talk given by the authors at the Mathematical Conference dedicated to Professor Ivan Prodanov on the
occasion of the 60${}^{th}$ anniversary of his birth and the 10${}^{th}$ anniversary of his death. The conference took place on May 16, 1995 at the Faculty of
Mathematics and Informatics of the Sofia University $``$St. Kliment Ohridski''.
It was planned  the talks of the participants in this conference to be published in a special volume of the Annuaire de l’Universite de Sofia $``$St. Kliment Ohridski”,
Faculte de Mathematiques et Informati\-que, Livre 1 - Mathematiques, but this has never happened. That's why we  have decided to publish our work separately. Since our files were lost and we had to write them once more,  the paper appears only now.

In 1979/80 academic year, Professor Ivan Prodanov organized a seminar on spectra at the Faculty of Mathematics of Sofia University. The participants in this
seminar, besides Iv. Prodanov, were G. Dimov, G. Gargov, Sv. Savchev, L. Stoyanov, V. Tchoukanov, T. Tinchev, D. Vakarelov. The talks of Iv. Prodanov on this
seminar were on his own investigations in the theory of abstract spectra and the uniqueness of Pontryagin-van Kampen duality. In the reviewing talks of the other
participants, Stone Duality Theorems for Boolean algebras and for distributive lattices (\cite{Sto1,Sto2}), H. A. Priestley's papers \cite{Pr1,Pr2,Pr3,Pr4}, M.
Hochster papers \cite{H1,H2}, the topological proof of G\"odel Completeness Theorem given by Rasiowa and Sikorski in \cite{RS} and many other interesting topics
were discussed.

Iv. Prodanov raised a number of interesting open problems on his seminar on spectra. Two of them were solved by some of the participants of the seminar and these
solutions caused, on their part, the appearence of other new papers. One of these problems was whether the category $\L_R$ of locally compact topological
$R$-modules, where $R$ is a locally compact commutative ring, admits precisely one (up to natural equivalence) functorial duality. (Using the classical
Pontryagin-van Kampen duality, one easily obtains a functorial duality in $\L_R$, called again Pontryagin duality. Hence, there is always a functorial duality in
$\L_R$.) L. Stoyanov \cite{St} showed that if $R$ is a compact commutative ring then the Pontryagin duality is the unique functorial duality in $\L_R$. Later on,
Gregorio \cite{Gre} and Gregorio and Orsatti \cite{GrOr94} generalized that result of Stoyanov. The second problem was whether a uniqueness theorem, like that for
Pontryagin-van Kampen duality, can be proved in the cases of Stone dualities for Boolean algebras and for distributive lattices. The answers were given by G. Dimov
in \cite{D1} and \cite{D2}, where it was proved that the Stone duality for Boolean algebras is unique and that there are only two (up to natural equivalence)
duality functors in the case of distributive lattices. Some very general results about representable dualities and the group of dualities were obtained  later on by G.
Dimov and W. Tholen in \cite{DT1,DT2}. It could be said that D. Vakarelov's paper \cite{V1} was also inspired by Prodanov's  seminar on spectra. This was certainly
so for the diploma thesis \cite{Sav1} of Sv. Savchev, written under the supervision of Professor Iv. Prodanov, and  for the paper \cite{Sav2}.

Iv. Prodanov presented his results on the uniqueness of Pontryagin-van Kampen duality in  the manuscripts \cite{P2} and \cite{P3}. The more than fifty-pages-long
paper \cite{P3} contained also an impressive list of open problems and conjectures. The publication of these manuscripts was postponed
 because Prodanov discovered that
analogous results were obtained earlier by D. Roeder \cite{R}. Prodanov's approach, however,  was different and even more general than that of D. Roeder.
 Only his untimely death withheld him from preparing
 these manuscripts for publication. The task of doing that was carried out by
  D. Dikranjan
and A. Orsatti. In their  paper \cite{DO}, all results of \cite{P2} and \cite{P3} were included and some of Prodanov's conjectures were answered. In such a way the
manuscripts \cite{P2} and \cite{P3} became known to the mathematical community and stimulated the appearance of other papers (see \cite{Dik94}, \cite{Greg89}).

The results of Iv. Prodanov on abstract spectra and separative algebras were announced in \cite{P1}, but their proofs were never written by him in the form of a
manuscript, preprint or paper. The very incomplete notes which we have from the Prodanov talks on the seminar on spectra seem to be the only trace of a small part of
these proofs. Since, in our opinion, the results, announced in \cite{P1}, are interesting and important, we decided to supply them with proofs. This is done in the
present paper, where we follow, in general, the exposition of \cite{P1}, but some of the announced there assertions are slightly generalized, some new statements
are added and some new applications are obtained. The main of the added results is Theorem \ref{th66} which was formulated and proved by us as a generalization of Prodanov's assertions Corollary \ref{th6} and Corollary \ref{th7}.

Section 1 of the paper is an introduction. Section 2, divided into four subsections, is devoted to the abstract spectra. In Subsection 2.1 the category \SS\ of
abstract spectra and their morphisms is introduced and studied. Subsection 2.2 contains two general examples of abstract spectra (see \ref{th2} and \ref{th4}). The
classical spectra of rings endowed with Zariski topology appear as special cases of the first of these examples (see \ref{ex1}), while the classical spectra of
distributive lattices with their Stone topology appear as special cases of both examples (see \ref{ex2} and \ref{ex20}). In Subsection 2.3 the main theorem of
Section 2 is proved (see \ref{th5}). This theorem asserts that the category \SS\ of abstract spectra and their morphisms is isomorphic to the category \CohSp\ of
coherent spaces and coherent maps and, hence, by the Stone Duality Theorem for distributive lattices, the category \SS\ is dual to the category \DLat\ of
distributive lattices and lattice homomorphisms. It is well known that the category \OStone\ of ordered Stone spaces and order-preserving continuous maps is also
dual to the category \DLat\ (see \cite{Pr1,Pr2} or \cite{J}), and that it is isomorphic to the category \CohSp\ (see, for example, \cite{J}).  Therefore, the
category \OStone\ is isomorphic to the category \SS. (The last fact could be also proved directly, but we don't do this.) So, each one of the categories \CohSp,
\OStone\ and \SS\ is dual to the category \DLat. In our opinion, the category \SS\ is the most natural and symmetrical one from all three of them. Subsection 2.4
contains two applications (see Corollary \ref{th6} and  Corollary \ref{th7}) of the already obtained results. The one from Corollary \ref{th6} is important for Section 3.
These applications appear  as special cases
of a general theorem (see Theorem \ref{th66}) which we formulate and prove here as a generalization of Prodanov's results Corollary \ref{th6} and  Corollary \ref{th7}.
Theorem \ref{th66} was used later on by us in our paper \cite{DV}.

   At the first glance the advent of spectra in so general   situations
as in \ref{th2} is unexpected, since psychologically they usually   are connected with separation. Actual\-ly, in general one does not know    whether   there are nontrivial prime
ideals but it turns out that   if the operations  $\times$  and  $+$   from   \ref{th2}   satisfy a   few   not   very restrictive natural conditions, then
the prime ideals become as many as in the commutative rings or in distributive lattices, for example.   In this way one comes to the notion of a {\em separative algebra}\/
considered in Section 3.

   Section 3 is divided into several   subsections.   In   Subsection
3.1  the   definition   of   a   {\it preseparative   algebra}\/  as   an algebra   with   two   multivalued   binary   opera\-tions  $\times$  and $+$
satisfying some natural axioms as commutativity and   associativity is given, and some calculus with this   operations is developed.   Subsection   3.2   is devoted to the theory
of   filters   and   ideals   in   preseparative algebras. The main notion of a  {\it separative   algebra}\/  is   given   in Subsection 3.3. Here a far of being
complete list of   examples   is given:      the   commutative     rings,     the distributive lattices   and   also   the   convex
spaces   (= separative algebras in which the two operations   coincide) are separative   algebras.    The main theorem for
separative algebras - the Separation theorem,   is proved in Subsection 3.4. In Subsection 3.5 some natural new   operations   in   separative algebras are studied and in Subsection 3.6 a general
representation   theorem for   separative   algebras   is   given.   Roughly     speaking,     every separative algebra $\underline{X}=(X, \times , +)$ can be embedded into a distributive
lattice $L$ in such a way that the operations in $\underline{X}$ are  obtained easily from the operations in $L$. That is new even   for   the   plane:   there   exists   a   distributive lattice
$L\supseteq \R^2$   such that for each segment $ab\subset \R^2$ one has $$ab=\{x\in\R^2 : x\le a\vee b\} = \{x\in\R^2 : x\ge a\wedge b\}.$$

The notion of separative   algebra   comes   from   an   analysis   of the separation theorems connected with   the   convexity.   The   abstract study   of   convexity
was started   by   Prenovitz   \cite{Pr}   and   different versions of the notion of {\em convex space}\/ appeared in \cite{P4},   \cite{P5}, \cite{T},   \cite{Br},   \cite{BrW}, \cite{PrJ}.
All they are   compared   in   \cite{Tes}.   The   convexity   was examined   from   other   aspects   in   \cite{B},    \cite{Can},  \cite{Gua}, \cite{Kay} and
\cite{Mah},   a   few applications are considered in  \cite{Zim}   and   \cite{BS}   contains   a critique.

   Y. Tagamlitzki \cite{T} obtained a general   Separation   theorem   for
convex spaces. It was   improved   (again   for   convex   spaces)   and applied to analytical separation problems in \cite{P4} and \cite{P5} (cf.   \cite{B} and
\cite{BrW}). It seems   however   that   the   natural   region   for   that theorem are not the convex spaces but   the   separative   algebras: the presence of
two operations makes the instrument more flexible, without additional   complications   (see   Subsection 3.4).     This permits to obtain   as special   cases   the   separation   by
prime   ideals   of   an ideal   and a multiplicative set in a commutative ring, or of an ideal and a filter in a distributive lattice, and also the separation of   two
convex sets by a convex set with convex complement.

   The paper ends with Subsection 3.8 devoted to a   generalization
of the Separation theorem for separative algebras supplied with a topology. Thus, even restricted   to   convex   spaces,   one   can   find, as in \cite{P5},   a   few
classical separation and representation theorems,               but   the presence of two operations  enlarges the possibilities for new applica\-tions.

Let's fix the notation. If $C$ denotes a category, we write $X\in |C|$ if $X$ is an object of $C$, and $f\in C(X,Y)$ if $f$ is a $C$-morphism with domain $X$ and
codomain $Y$. All lattices will be with top (=unit) and bottom (=zero) elements, denoted respectively by 1 and 0. We don't require the elements $0$ and $1$ to be
distinct. As usual, the lattice homomorphisms are assumed to preserve the distinguished elements $0$ and $1$. \DLat\ will stand for the category of  distributive
lattices and lattice homomorphisms. If $X$ is a set then we write $Exp(X)$ for the set of all subsets of $X$ and denote by $|X|$ the cardinality of $X$. If $(X,\T)$
is a topological space and $A$ is a subset of $X$ then $cl_{(X,\T)}A$ or, simply, $cl_XA$ stands for the closure of $A$ in the space $(X,\T)$. We denote by \D\ the
two-point discrete topological space and by \Set\ the category of all sets and functions between them. As usual, we say that a preordered set $(X,\le)$ (i.e. $\le$
is a reflexive and transitive binary relation on $X$) is a {\it directed set} (resp. an {\it ordered set}) if for any $x,y\in X$ there exists a $z\in X$ such that
$x\le z$ and $y\le z$ (resp. if the relation $\le$ is also antisymmetric).

Our main references are: \cite{J} -- for category theory and Stone dualities,  \cite{E} --  for general topology,  and \cite{L} -- for algebra.


\section{Spectra}


\subsection{The Category of Abstract Spectra}

\begin{nota}\label{no1} \rm Let $(S,\T^+,\T^-)$ be a non-empty bitopological space. Then we put $\L^+ = \{U\in \T^+ : S\stm U\in\T^-\}$ and $\L^- = \{U\in \T^- :
S\stm U\in\T^+\}$. \end{nota}

\begin{pro}\label{pr0} Let $(S,\T^+,\T^-)$ be a non-empty bitopological space. Then the families $\L^+$ and $\L^-$ (see \ref{no1} for the notation) are closed
under finite unions and finite intersections. \end{pro}

\doc It is obvious. \sq

\begin{defi}\label{de1} \rm A non-empty bitopological space $(S,\T^+,\T^-)$ is called an {\em abstract spectrum}, if it has the following properties:

\hspace{-1cm}(SP1) $\L^+$ is a base for $\T^+$ and $\L^-$ is a base for $\T^-$;

\hspace{-1cm}(SP2) if $F\sbe S$ and $S\stm F\in\T^+$ (resp. $S\stm F\in\T^-$) then $F$ is a compact subset of the topological space $(S,\T^-)$ (resp. $(S,\T^+)$);

\hspace{-1cm}(SP3) at least one of the topological spaces $(S,\T^+)$ and $(S,\T^-)$ is a $T_0$-space. \end{defi}

\begin{pro}\label{pr1} If $(S,\T^+,\T^-)$ is an abstract spectrum then $(S,\T^+)$ and $(S,\T^-)$ are compact $T_0$-spaces. \end{pro}

\doc By (SP3), one of the spaces $(S,\T^+)$ and $(S,\T^-)$ is $T_0$-space. Let, for example, $(S,\T^+)$ be  a $T_0$-space. Then we shall prove that $(S,\T^-)$ is
also a $T_0$-space.

Let $x,y\in S$ and $x\not= y$. Then there exists $U\in \T^+$ such that $|U\cap\{x,y\}|=1$. Let, for example, $x\in U$. Then, using (SP1), we can find a $V\in\L^+$
such that $x\in V\sbe U$. Putting $W=S\stm V$, we obtain that $W\in \T^-$, $y\in W$ and $x\not\in W$. Therefore, $(S,\T^-)$ is a $T_0$-space.

Since $S$ is a closed subset of $(S,\T^+)$, the condition (SP2) implies that $S$ is a compact subset of $(S,\T^-)$.

Analogously, we obtain that $(S,\T^+)$ is a compact space. \sq

\begin{pro}\label{pr2} Let $(S,\T^+,\T^-)$ be an abstract spectrum. Then $\L^+ = \{U\in \T^+: U$   is a compact subset of  $(S,\T^+)\}$ and $\L^- = \{U\in \T^-: U$
is   a compact   subset
  of  $(S,\T^-)\}$
(see \ref{no1} for the notation). \end{pro}

\doc Let us prove first that $\L^+ = \{U\in \T^+: U$   is a compact subset of  $(S,\T^+)\}$.

If $V\in\L^+$ then $S\stm V\in\T^-$. Hence $V$ is a closed subset of $(S,\T^-)$.            This implies, by (SP2), that $V$ is a compact subset of  $(S,\T^+)$.
Conversely, if $U\in\T^+$ and $U$ is a compact subset of  $(S,\T^+)$ then for every $x\in U$ there exists a $U_x\in\L^+$ such that $x\in U_x\sbe U$. Choose a finite
subcover $\{U_{x_i}:\ion\}$ of the cover $\{U_x:x\in U\}$ of the compact set $U$. Then $U= \bigcup\{U_{x_i}:\ion\}$ and hence, by \ref{pr0}, $U\in\L^+$.

The proof of the equation $\L^- = \{U\in \T^-: U$   is   a compact   subset
  of  $(S,\T^-)\}$
is analogous. \sq

\begin{pro}\label{pr3} Let $(S,\T^+,\T^-)$ be an abstract spectrum. Then $\L^+ = \{S\stm U: U\in\L^-\}$ and $\L^- = \{S\stm U: U\in\L^+\}$ (see \ref{no1} for the
notation). \end{pro}

\doc Let's prove that $\L^- = \{S\stm U: U\in\L^+\}$.

Take $V\in \L^-$ and put $U=S\stm V$. Then $U\in\T^+$ and $S\stm U\in\L^- \sbe\T^-$. Hence, $U\in\L^+$ and $V=S\stm U$. Conversely, if $U\in\L^+$ then $V=S\stm
U\in\T^-$ and $S\stm V\in\L^+\sbe\T^+$. Therefore, $S\stm U\in \L^-$.

The proof of the equation $\L^+ = \{S\stm U: U\in\L^-\}$  is analogous. \sq

\begin{cor}\label{co1} Let $(S,\T^+,\T^-_1)$ and $(S,\T^+,\T^-_2)$  be abstract spectra. Then the topolo\-gies $\T^-_1$ and $\T^-_2$ coincide. \end{cor}

\doc It follows directly from \ref{pr2}, \ref{pr3} and (SP1) (see \ref{de1}). \sq

\begin{defi}\label{de2} \rm Let $(S_1,\T^+_1,\T^-_1)$ and $(S_2,\T^+_2,\T^-_2)$  be abstract spectra. Then a function $f\in\Set(S_1,S_2)$ is called an
$\SS$\/{\em-morphism} if $f: (S_1,\T^+_1)\lra (S_2,\T^+_2)$ and $f: (S_1,\T^-_1)\lra (S_2,\T^-_2)$ are continuous maps. The class of all abstract spectra together
with the class of all $\SS$\/-morphisms and the natural composition between them form, obviously, a category which will be denoted by $\SS$ and will be called {\em
the category of abstract spectra}. \end{defi}

\begin{defi}\label{de3} \rm An abstract spectrum $(S,\T^+,\T^-)$ is called a {\em Stone spectrum} if the topologies $\T^+$ and $\T^-$ coincide. \end{defi}

\begin{pro}\label{pr4} Let $(S,\T)$ be a topological space. Then the bitopological space $(S,\T,\T)$ is a Stone spectrum iff $(S,\T)$ is a Stone space. \end{pro}

\doc $(\Ra)$ Let $(S,\T,\T)$ be a Stone spectrum. Then, by \ref{pr1}, $(S,\T)$ is a compact $T_0$-space. According to (SP1) (see \ref{de1}), the family
$\L^+=\{U\in\T:S\stm U\in \T\}$ is a base for $\T$. Consequently $(S,\T)$ is a zero-dimensional space. We shall show that it is also a $T_2$-space. Indeed, let
$x,y\in S$ and $x\not=y$. Then there exists a $U\in\T$ such that $|U\cap\{x,y\}|=1$. Let, for example, $x\in U$. Since $\L^+$ is a base for $\T$, we can find a
$V\in\L^+$ such that $x\in V\sbe U$. Then $x\in V\in\T$ and $y\in S\stm V\in\T$. Therefore, $(S,\T)$ is a $T_2$-space. So, we proved that $(S,\T)$ is a compact
zero-dimensional $T_2$-space, i.e. a Stone space.

 \smallskip

 \noindent$(\La)$ Let $(S,\T)$ be a Stone space. Put $\L=\{U\in\T:S\stm U\in\T\}$ and $\T^+=\T^-=\T$. Then $\L^+=\L=\L^-$ (see \ref{no1} for the notation). We shall prove
that $(S,\T^+,\T^-)$ is an abstract spectrum. Then it will be automatically a Stone spectrum. Since $\L$ is a base for $(S,\T)$, the axiom (SP1) (see \ref{de1}) is
fulfilled. The axioms (SP2) and (SP3) are also fulfilled, since $(S,\T)$ is a compact $T_2$-space. Consequently $(S,\T^+,\T^-)$ is an abstract spectrum. \sq

\begin{pro}\label{pr5} An abstract spectrum $(S,\T^+,\T^-)$ is a Stone spectrum iff $(S,\T^+)$ and $(S,\T^-)$ are $T_1$\/-spaces. \end{pro}

\doc $(\Ra)$ Since $(S,\T^+,\T^-)$ is a Stone spectrum, we have that $\T^+=\T^-$. Then \ref{pr4} implies that $(S,\T^+)$ and $(S,\T^-)$      are even $T_2$-spaces.

 \smallskip

 \noindent$(\La)$  Let $(S,\T^+)$     and $(S,\T^-)$     are $T_1$-spaces. We shall prove that $\T^+=\T^-$.

Let $U\in\T^-$. Then $S\stm U$ is closed in $(S,\T^-)$     and hence, by \ref{pr1}, it is a compact subset of $(S,\T^-)$. Let $x\in U$. Since $(S,\T^+)$     is a
$T_1$-space, for every $y\in S\stm U$ there exists a $V_y\in\L^+$ such that $x\in V_y\sbe S\stm\{y\}$. Hence $y\in S\stm V_y\sbe S\stm\{x\}$ and $S\stm V_y\in\T^-$.
Let $\{S\stm V_{y_i}:\ion\}$ be a finite subcover of the cover $\{S\stm V_y:y\in S\stm U\}$ of $S\stm U$ and let $V_x=\bigcap\{V_{y_i}:\ion\}$. Then $x\in V_x\in
\T^+$ and $V_x\sbe U$. We obtain that $U=\bigcup\{V_x:x\in U\}\in\T^+$. Hence $\T^-\sbe\T^+$. Analogously, using the fact that $(S,\T^-)$     is a $T_1$-space, we
prove that $\T^+\sbe \T^-$. Therefore $\T^+=\T^-$, i.e. $(S,\T^-,\T^+)$  is a Stone spectrum. \sq

\begin{rem}\label{re1} \rm Let $(S,\T^+,\T^-)$ be an abstract spectrum. Then, arguing as in \ref{pr1}, we obtain that $(S,\T^+)$ and $(S,\T^-)$ are $T_1$\/-spaces
iff at least one of them is a $T_1$\/-space. \end{rem}

\begin{pro}\label{pr6} Let $(S,\T^+,\T^-)$ be an abstract spectrum and let us put $\T = \sup\{\T^+,\T^-\}$. Then $(S,\T)$ is a Stone space and hence (see \ref{pr4})
$(S,\T,\T)$ is a Stone spectrum. \end{pro}

\doc The topology $\T$ has as a subbase the family $\PP=\T^+\cup\T^-$. Hence the family $\BB=\{U^+\cap U^-:U^+\in \T^+,\ U^-\in\T^-\}$ is a base for $\T$. Then,
obviously, the family $\BB_0=\{U^+\cap U^-:U^+\in \L^+,\ U^-\in\L^-\}$ is also a base for $\T$. For every $U\in\L^+$ we have that $U\in\T^+\sbe\T$ and $S\stm
U\in\T^-\sbe \T$. Consequently the elements of $\L^+$ are clopen subsets of $(S,\T)$. Obviously, the same is true for the elements of $\L^-$. Hence the elements of
$\BB_0$ are clopen in $(S,\T)$, which implies that $(S,\T)$ is a zero-dimensional space. This fact, together with (SP3) (see \ref{de1}), shows that $(S,\T)$ is a
Hausdorff space.

Applying Alexander subbase theorem to the subbase $\PP$ of $(S,\T)$, we shall prove that $(S,\T)$ is a compact space. Indeed, let $S=\bigcup\{U_\a\in\T^+:\a\in A\}\
\cup\ \bigcup\{V_\b\in\T^-:\b\in B\}$ and $F=S\stm\bigcup\{U_\a:\a\in A\}$. Then $F\sbe\bigcup\{V_\b:\b\in B\}$ and $F$ is closed in $(S,\T^+)$. Consequently, by
(SP2) (see \ref{de1}), $F$ is a compact subset of $(S,\T^-)$. This implies that there exist $\b_1,\ldots,\b_n\in B$ such that $F\sbe\bigcup\{V_{\b_i}:\ion\}$. Then
$G=S\stm\bigcup\{V_{\b_i}:\ion\}\sbe\bigcup\{U_{\a}:\a\in A\}$. Since $G$ is a closed subset of $(S,\T^-)$, it is a compact subset of $(S,\T^+)$ (by (SP2)). Hence,
there exist $\a_1,\ldots,\a_m\in A$ such that $G\sbe\bigcup\{U_{\a_j}:\jom\}$. Therefore, $S=\bigcup\{U_{\a_j}:\jom\}\ \cup\ \bigcup\{V_{\b_i}:\ion\}$.  This shows
that $(S,\T)$ is compact. Hence, $(S,\T)$ is a Stone space. \sq

\begin{rem}\label{re2} \rm Let $(S,\T^+,\T^-)$ be an abstract spectrum and $id:S\lra S$, $x\lra x$, be the identity function. Then, obviously,
$id\in\SS((S,\T,\T),(S,\T^+,\T^-))$ (see \ref{pr6} for the notation). \end{rem}

\begin{pro}\label{pr61} Let $(S,\T^+,\T^-)$ be a bitopological space such that $\L^+$ is a base for $\T^+$ and $\L^-$ is a base for $\T^-$ (see \ref{no1} for the
notation). Let $\T = \sup\{\T^+,\T^-\}$, $(S,\T)$ be a compact $T_2$-space, $S_1\sbe S$, $\T^+_1 = \{U\cap S_1:U\in\T^+\}$ and $\T^-_1 = \{U\cap S_1:U\in\T^-\}$.
Then the bitopological space $(S_1,\T^+_1,\T^-_1)$ is an abstract spectrum iff $S_1$ is a closed subset of the topological space $(S,\T)$. \end{pro}

\doc  $(\Ra)$ Let $\T_1 = \sup\{\T^+_1,\T^-_1\}$. Then, by \ref{pr6}, $(S_1,\T_1)$ is a Stone space. Hence it is a compact Hausdorff space. Since, obviously,
$\T_1=\T |S_1$, we obtain that $S_1$ is a compact subspace of the Hausdorff space $(S,\T)$. Consequently $S_1$ is a closed subset of $(S,\T)$.

 \smallskip

 \noindent$(\La)$ We shall show that $(S_1,\T^+_1,\T^-_1)$ is an abstract spectrum. Let $\L^+_1 = \{U\cap S_1:U\in\L^+\}$, $\L^-_1 = \{U\cap S_1:U\in\L^-\}$, $\L^+_{S_1} =
\{U\in \T^+_1 : S_1\stm U\in\T^-_1\}$ and $\L^-_{S_1} = \{U\in \T^-_1 : S_1\stm U\in\T^+_1\}$. Then, obviously, $\L^+_1\sbe\L^+_{S_1}$ and $\L^-_1\sbe\L^-_{S_1}$.
Since $\L^+_1$  (resp. $\L^-_1$) is a base for $(S_1,\T^+_1)$ (resp. $(S_1,\T^-_1)$), we obtain that $\L^+_{S_1}$ (resp. $\L^-_{S_1}$) is a base for $(S_1,\T^+_1)$
(resp. $(S_1,\T^-_1)$). Hence  the condition (SP1) (see \ref{de1}) is fulfilled.

In the part $(\Ra)$ of this proof, we  noted that the topology $\T_1 = \sup\{\T^+_1,\T^-_1\}$ on $S_1$ coincides with the topology $\T |S_1$. Hence, $(S_1,\T_1)$ is
a compact Hausdorff space (since $(S,\T)$ is such and $S_1$ is a closed subset of $(S,\T)$). Let now $F$ be a closed subset of $(S_1,\T^+_1)$ (resp.
$(S_1,\T^-_1)$). Then $F$ is a closed subset of $(S_1,\T_1)$. Therefore $F$ is a compact subset of $(S_1,\T_1)$. Since the identity maps $id: (S_1,\T_1)\lra
(S_1,\T_1^+)$ and $id: (S_1,\T_1)\lra (S_1,\T_1^-)$ are continuous, we obtain that $F$ is a compact subset of $(S_1,\T_1^-)$ (resp. $(S_1,\T_1^+)$). Hence, the
condition (SP2) (see \ref{de1}) is fulfilled.

For showing that the condition (SP3) (see \ref{de1}) is fulfilled, it is enough to prove that $(S_1,\T^+_1)$ is a $T_0$-space. Let $x,y\in S_1$ and $x\not=y$. Since
$(S_1,\T_1)$ is a $T_2$-space, there exist $U\in\L^+_1$ and $V\in\L^-_1$ such that $x\in U\cap V\sbe S_1\stm\{y\}$. If $y\not\in U$ then  the element $U$ of
$\T^+_1$ separates $x$ and $y$. If $y\in U$ then $y\not\in V$. Hence $y\in S_1\stm V$ and $x\not\in S_1\stm V$. Since $S_1\stm V\in\T^+_1$, we obtain that $x$ and
$y$ are separated by an element of $\T^+_1$. Consequently, $(S_1,\T^+_1)$ is a $T_0$-space. \sq

\begin{cor}\label{co2} Let $(S,\T^+,\T^-)$ be an abstract spectrum, $\T = \sup\{\T^+,\T^-\}$, $S_1\sbe S$, $\T^+_1 = \{U\cap S_1:U\in\T^+\}$ and $\T^-_1 = \{U\cap
S_1:U\in\T^-\}$. Then the bitopological space $(S_1,\T^+_1,\T^-_1)$ is an abstract spectrum iff $S_1$ is a closed subset of the topological space $(S,\T)$.
\end{cor}

\doc It follows immediately from \ref{pr61}, \ref{de1} and \ref{pr6}. \sq


\subsection{Examples of Abstract Spectra}

 \begin{lm}\label{lm20}
 Let $X$ be a set and $Exp(X)$ be the family of all subsets of $X$.
 Let's put, for every $x\in X$,
 $\tilde U^+_x = \{A\sbe X: x\not\in A\}$ and
 $\tilde U^-_x = \{A\sbe X: x\in A\}$.
 Let $\tilde \PP^+ = \{\tilde U^+_x : x\in X\}$,
 $\tilde\PP^- = \{\tilde U^-_x : x\in X\}$,
 $\tilde\T^+$ (resp. $\tilde\T^-$) be the topology on $Exp(X)$
 having $\tilde\PP^+$ (resp. $\tilde\PP^-$)
 as a subbase and
 $\tilde\T = \sup\{\tilde\T^+,\tilde\T^-\}$. Let's identify the set $Exp(X)$
 with the set $\D^X$
 (where \D\ is the two-point set $\{0,1\}$)
 by means of the map $e:Exp(X)\lra \D^X$, $A\sbe X\lra \chi_A$, where
 $\chi_A:X\lra\D$ is the characteristic function of $A$, i.e.
 $\chi_A(x)=1$ if $x\in A$ and
 $\chi_A(x)=0$ if $x\not\in A$. Then
 the topology
 $\tilde\T$ on $Exp(X)$ coincides with the Tychonoff topology on $\D^X$
 (where the set $\D$ is endowed with the discrete topology).
 \end{lm}

 \doc Let $\tilde \PP=\tilde \PP^+\cup\tilde \PP^-$. Then
 $\tilde \PP$ is a subbase for the topology $\tilde \T$ on $Exp(X)$.
 For every $x\in X$ we have, identifying $Exp(X)$ and $\D^X$ by means of the
 map $e$, that
 $\tilde U^+_x=\{f\in\D^X:f(x)=0\}$ and
 $\tilde U^-_x=\{f\in\D^X:f(x)=1\}$. Now it becomes clear that
 the family $\tilde \PP$ is also a subbase for the Tychonoff topology on
 $\D^X$ when $\D$ is endowed with the discrete topology. Therefore
 the topology
 $\tilde\T$ on $Exp(X)$ coincides with the Tychonoff topology on $\D^X$.
\sq

 \begin{pro}\label{th1}
 Let $X$ be a set and $S$ be family of subsets of $X$ (i.e. $S\sbe Exp(X)$).
 Let's put, for every $x\in X$,
 $U^+_x = \{p\in S: x\not\in p\}$ and
 $U^-_x = \{p\in S: x\in p\}$.
 Let $\PP^+ = \{U^+_x : x\in X\}$,
 $\PP^- = \{U^-_x : x\in X\}$,
 $\T^+$ (resp. $\T^-$) be the topology on $S$ having $\PP^+$ (resp. $\PP^-$)
 as a subbase and
 $\T = \sup\{\T^+,\T^-\}$.
 Then the following conditions are equivalent:

 \smallskip

 \noindent(a) $(S,\T^+,\T^-)$ is an abstract spectrum;

 \smallskip

 \noindent(b) $(S,\T)$ is a compact $T_2$-space;

 \smallskip

 \noindent(c) $S$ is a closed subset
 of the Cantor cube $\D^X$
 (where \D\ is the discrete two-point space and $S$ is identified with a
 subset of $\D^X$ as in \ref{lm20}).
 \end{pro}

 \doc $(a)\Ra (b).$ This follows from \ref{pr6}.

 \smallskip

 \noindent$(b)\Ra (a).$ Let $x\in X$. Then $S\stm U^+_x=U^-_x$ and $S\stm U^-_x=U^+_x$.
 Hence $\PP^+\sbe \L^+$ and $\PP^-\sbe \L^-$ (see \ref{no1} for the
 notation). Consequently, using \ref{pr0}, we obtain that $\L^+$ (resp. $\L^-$)
 is a base for $(S,\T^+)$ (resp. $(S,\T^-)$). This shows that putting
 $S_1=S$ in \ref{pr61}, we get that $(S,\T^+,\T^-)$ is an abstract
 spectrum.

 \smallskip

 \noindent$(b)\Ra(c).$ It is clear from the corresponding definitions that, using the
 notation of \ref{lm20}, we have $\tilde U^+_x\cap S=U^+_x$ and
 $\tilde U^-_x\cap S=U^-_x$ for every $x\in X$. Hence, by \ref{lm20},
 the topology $\T$ on $S$ coincides with the subspace topology on $S$
 induced by the Tychonoff topology on $\D^X$.
 Then the condition (b) and the
 fact that $\D^X$ is a Hausdorff space imply that $S$ is a closed subset
 of the Cantor cube $\D^X$.

 \smallskip

 \noindent$(c)\Ra(b).$ In the preceding paragraph we have already noted that
 the topology $\T$ on $S$ coincides with the subspace topology on $S$
 induced by the Tychonoff topology on $\D^X$.
 Therefore the condition (c) implies that $(S,\T)$ is a compact Hausdorff
 space (since $\D^X$ is such).
 \sq

 \begin{defi}\label{de4}
 \rm
 Let $X$ be a set endowed with two arbitrary multivalued binary operations
 $\op$ and $\ot$.
 Let's call a subset $p$ of $X$ a {\em prime ideal in}
 $(X,\op,\ot)$ if the following two conditions are fulfilled:

 \smallskip

 \noindent i) if $x,y\in p$ then $x\op y\sbe p$;

 \smallskip

 \noindent ii) if $(x\ot y)\cap p\not =\ems$ then $x\in p$ or $y\in p$.

 Let us fix two
 different points
 $0$ and $1$
 of $X$.
 We shall say that
 a prime ideal $p\sbe X$ is {\em proper}
 (or, more precisely, {\em proper with respect to the points $0$ and $1$}),
 if $0\in p$ and $1\not\in p$.

 A subset $q$ of $X$ is called a {\em prime (proper) flter in}
 $(X,\op,\ot)$ if the set $X\stm q$ is a prime (proper) ideal.
 \end{defi}

 \begin{theorem}\label{th2}
 Let $X$ be a set endowed with two arbitrary multivalued binary operations
 $\op$ and $\ot$
 and two fixed different points $\xi_o$ and $\xi_1$.
 Denote by $S(X)$ (resp. $S(X)_{pr}$) the set of all (resp. all
 proper) prime ideals in
 $(X,\op,\ot)$ and define the topologies $\T^+$ and $\T^-$ on $S(X)$
 (resp. $\T^+_{pr}$ and $\T^-_{pr}$ on $S(X)_{pr}$)
 exactly as in \ref{th1}. Then the bitopological spaces
 $(S(X),\T^+,\T^-)$ and $(S(X)_{pr},\T^+_{pr},\T^-_{pr})$ are abstract spectra.
 \end{theorem}

 \doc We  first prove  that
 the bitopological space $(S(X),\T^+,\T^-)$ is an abstract spectrum.
 For doing this it is enough to show that
 $S(X)$ is a closed subset of the Cantor cube $\D^X$ (see \ref{th1}).

 Let $\{p_\s\in S(X):\s\in \Sigma\}$ be a net in the Cantor cube $\D^X$
 converging to a point $p\in\D^X$. We have to prove that $p\in S(X)$, i.e.
 that $p$ is a prime ideal in $(X,\op,\ot)$.
 Let $f_\s=e(p_\s)$ and
 $f=e(p)$ (see \ref{lm20} for the notation). Then the net
 $\{f_\s,\s\in\Sigma\}$ in $\D^X$ converges to $f$, i.e., for every $x\in X$,
 the net $\{f_\s(x),\s\in\Sigma\}$ in the discrete space $\D$ converges to
 $f(x)$.

 Let $a,b\in p$. Then $f(a)=f(b)=1$. Therefore there exists
 a $\s_0\in\Sigma$ such that $f_\s(a)=1=f_\s(b)$ for every $\s>\s_0$.
 This means that for every $\s>\s_0$ we have that $a\in p_\s$ and
 $b\in p_\s$. Since $p_\s$ is a prime ideal, we obtain that
 $a\op b\sbe p_\s$ for every $\s>\s_0$. Then, for every $x\in a\op b$ and
 for every $\s>\s_0$, we have that $f_\s(x)=1$. This implies that
 $f(x)=1$ for every $x\in a\op b$. Hence,
 if $x\in a\op b$ then
 $x\in p$,
 i.e.
 $a\op b\sbe p$.

 Let $a,b\in X$ and $(a\ot b)\cap p\not=\ems$. Then there exists
 a $x\in (a\ot b)\cap p$. Hence $f(x)=1$. This implies that there exists
 a $\s_0\in\Sigma$ such that $f_\s(x)=1$ for every $\s>\s_0$.
 Consequently $x\in p_\s$ for every $\s>\s_0$. Then $(a\ot b)\cap p_\s\not=\ems$
 for every $\s>\s_0$. Hence, for every $\s>\s_0$, we have that
 $a\in p_\s$ or $b\in p_\s$, i.e.   $f_\s(a)=1$ or
 $f_\s(b)=1$. Suppose that $a\not\in p$ and $b\not\in p$. Then $f(a)=0=f(b)$.
 Therefore, there exists a $\s_1\in\Sigma$ such that $f_\s(a)=f_\s(b)=0$
 for every $\s>\s_1$. Since for every $\s>\sup\{\s_0,\s_1\}$ we have
 that   $f_\s(a)=1$ or $f_\s(b)=1$, we get a contradiction.
 Hence we obtain that $a\in p$ or $b\in p$. So, we proved that $p$ is a
 prime ideal in $(X,\op,\ot)$.
 This shows that $S(X)$ is a closed subset
 of the Cantor cube $\D^X$. Hence,
 the bitopological space $(S(X),\T^+,\T^-)$ is an abstract spectrum.

 If the prime ideals $p_\s$ in the above proof were proper, then,
 obviously, $p$ would be  also  proper.
 This shows that the set $S(X)_{pr}$ is also
 a closed subset of the Cantor cube $\D^X$. So, the bitopological space
 $(S(X)_{pr},\T^+_{pr},\T^-_{pr})$ is an abstract spectrum.
 \sq

 \begin{exa}\label{ex1}
 \rm
 Let $(A,+,.)$ be a commutative ring with unit $(0\not = 1)$,
 $x\op y$ be the ideal in the ring $(A,+,.)$ generated by $\{x,y\}$, and
 $x\ot y = x.y$, for every $x,y\in A$.
 Then, applying the construction
 from \ref{th2} to the set $A$ with the operations $\op$ and $\ot$
 and with fixed points $0$ and $1$, we get
 the topological space $(S(A)_{pr},\T^+_{pr})$. We assert that it coincides
 with the classical spectrum of the ring $(A,+,.)$.
 \end{exa}

 \doc Recall that: a) a subgroup
 $I$ of the additive group $(A,+)$ is called an
 {\em ideal} in the
 commutative ring $(A,+,.)$ with unit if $A.I=I$; b) an ideal
 $p\not = A$ in the ring $A$ is said to be a {\em prime
 ideal} if $(x,y\in A,\ x.y\in p)\Ra$
 $(x\in p$ or $y\in p)$; c) the set of all
 prime ideals in the commutative ring $A$ is denoted by $spec(A)$; d) the
 family ${\cal Z} =
 \{U_I=\{p\in spec(A): I\not\sbe p\}: I$ is an ideal in $A\}$
 is a topology on the set $spec(A)$, called {\em Zariski topology};
 e) the topological
 space $(spec(A), {\cal Z})$ is the classical spectrum of the commutative ring
 $(A,+,.)$ with unit.

 We shall denote by $I(M)$ the ideal in $A$ generated
 by a subset $M$ of  $A$.

 We first prove that the sets $spec(A)$ and $S(A)_{pr}$ coincide.

 Let $p\in S(A)_{pr}$. Then $1\not\in p$ and hence $p\not = A$.
 If $a,b\in p$ then $a\op b\sbe p$, i.e. $I(\{a,b\})\sbe p$. Hence, $a-b\in p$.
 This shows that $p$ is an additive subgroup of $A$.
 Let $x\in A$ and $a\in p$.
 Since $a\op a = I(\{a\})\sbe p$, we get that $x.a\in p$.
 If $x,y\in A$ and $x.y\in p$, then $(x\ot y)\cap p\not =\ems$ and, hence,
   $x\in p$ or $y\in p$. Consequently, we proved that $p\in spec(A)$.

 Conversely, let $p\in spec(A)$ and $a,b\in p$. Then, obviously,
 $I(\{a,b\})\sbe p$ and, hence, $a\op b\sbe p$. If $(a\ot b)\cap p\not=\ems$
 then $a.b\in p$.
 This implies that $a\in p$ or $b\in p$. Since $1\not\in p$, we get
 that $p\in S(A)_{pr}$. Therefore, $S(A)_{pr}=spec(A)$.

 Now we prove that $\T^+_{pr}=\cal Z$.

 Let $a\in A$. Then, obviously, $U^+_a=\{p\in S(A)_{pr}: a\not\in p\}=
 \{p\in spec(A):I(\{a\})\not\sbe p\}\in\cal Z$. Hence, $\T^+_{pr}\sbe\cal Z$.
 Conversely, let $U\in\cal Z$. Then there exists an ideal $I$ in $A$ such
 that $U=U_I$. Let $p\in U$. Then there exists an
 $a=a(p)\in I\stm p$. Hence
 $p\in U^+_a$. We shall prove that $U^+_a\sbe U$. Indeed, if $q\in U^+_a$
 then $a\not\in q$ and, consequently,
 $I\not\sbe q$. This shows that $q\in U_I=U$. So, we obtained that
 $p\in U^+_a\sbe U$. Therefore, ${\cal Z}\sbe \T^+_{pr}$.
 \sq

 \begin{exa}\label{ex2}
 \rm
 Let $(L,\vee,\we)$ be a distributive lattice with $0$ and $1$
 and let us put
 $x\op y=\{z\in L: z\le x\vee y\}$ and
 $x\ot y=\{z\in L: z\ge x\we y\}$, for every $x,y\in L$.
 Then, applying the construction
 from \ref{th2} to the set $L$ with the operations $\op$ and $\ot$
 and with fixed points $0$ and $1$, we get
 the topological space $(S(L)_{pr},\T^+_{pr})$. We assert that it coincides
 with the classical spectrum $spec(L)$
 of the distributive lattice $(L,\vee,\we)$.
 \end{exa}

 \doc Recall that: a) a sub-join-semillatice $I$ of the lattice $L$ is said
 to be an {\em ideal}
 in $L$ if $(a\in I,\ b\in L$ and $b\le a)\Ra (b\in I)$; b)
 an ideal $p$ in $L$ is called a {\em prime ideal} if $1\not\in p$ and
 $(a\we b\in p)\Ra$   $(a\in p$ or $b\in p)$; c) the set of all
 prime ideals in $L$ is denoted by $spec(L)$; d)
 the family
 ${\cal O} = \{U_I=\{p\in spec(L): I\not\sbe p\}: I$ is an ideal in $L\}$
 is a topology on the set $spec(L)$, called {\em Stone topology};
 e) the topological
 space $(spec(L), {\cal O})$ is the classical spectrum of the lattice
 $(L,\vee,\we,0,1)$.

 We  first prove that the sets $spec(L)$ and $S(L)_{pr}$ coincide.

 Let $p\in S(L)_{pr}$. Then $0\in p$ and $1\not\in p$. If $a,b\in p$ then
 $a\op b\sbe p$ and, hence, $a\vee b\in p$. Let $c\in L$,
 $a\in p$ and $c\le a$.
 Since $a\in p$, we have that $a\op a\sbe p$ and, consequently, $c\in p$. If
 $c,d\in L$ and $c\we d\in p$ then $(c\ot d)\cap p\not=\ems$. Therefore
 $c\in p$ or $d\in p$. So, $p\in spec(L)$.

 Let $p\in spec(L)$ and $a,b\in p$. Then $a\vee b\in p$ and, for all $c\in L$
 such that $c\le a\vee b$, we have that $c\in p$. Hence $a\op b\sbe p$.
 Let $x,y\in p$ and $(x\ot y)\cap p\not =\ems$. Then there exists a $z\in p$
 such that $z\ge x\we y$. Hence $x\we y\in p$. This implies that
 $x\in p$ or $y\in p$. Since $1\not\in p$, we obtain that $p\in S(L)_{pr}$.
 So, $S(L)_{pr}=spec(L)$.

 Now we prove that $\T^+_{pr}=\cal O$.

 Let $a\in L$ and $I(a)=\{x\in L:x\le a\}$. Then  $I(a)$ is
 an ideal in $L$.  Obviously,
 $U^+_a=\{p\in S(L)_{pr}: a\not\in p\}=\{p\in spec(L): I(a)\not\sbe p\}\in\cal O$.
 Hence $\T^+_{pr}\sbe\cal O$. Conversely, let $U\in \cal O$. Then there exists
 an ideal $I$ in $L$ such that $U=U_I$.
 Let $p\in U$. Then there exists
 an $a=a(p)\in I\stm p$. Hence $p\in U^+_a$ and we need to prove only that
 $U^+_a\sbe U$. Let $q\in U^+_a$. Then $a\not\in q$. Consequently $I\not\sbe q$,
 which means that $q\in U_I=U$. So, $p\in U^+_a\sbe U$. We obtained that
 ${\cal O}\sbe \T^+_{pr}$.
 \sq

 \begin{defi}\label{de5}
 \rm
 Let $X$ be a set endowed with two arbitrary single-valued binary operations
 $+$ and $\times$.
 Let's call a subset $p$ of $X$ an {\em l-prime ideal in}
 $(X,+,\times)$ if the following two conditions are fulfilled:

 \smallskip

 \noindent i) $x+y\in p$ iff $x\in p$ and $y\in p$;

 \smallskip

 \noindent ii) $x\times y\in p$ iff   $x\in p$ or $y\in p$.

 Let us fix two
 different points
 $0$ and $1$
 of $X$.
 We shall say that
 an l-prime ideal $p\sbe X$ is {\em proper}
 (or, more precisely, {\em proper with respect to the points $0$ and $1$}),
 if $0\in p$ and $1\not\in p$.
 \end{defi}

 \begin{theorem}\label{th4}
 Let $X$ be a set endowed with two arbitrary single-valued binary operations
 $+$ and $\times$
 and two fixed different points $\xi_o\in X$ and $\xi_1\in X$.
 Denote by $S'(X)$ (resp. $S'(X)_{pr}$) the set of all (proper) l-prime
 ideals in
 $(X,+,\times)$ and define the topologies $\T^+$ and $\T^-$ on $S'(X)$
 (resp. $\T^+_{pr}$ and $\T^-_{pr}$ on $S'(X)_{pr}$)
 exactly as in \ref{th1}. Then the bitopological spaces
 $(S'(X),\T^+,\T^-)$ and $(S'(X)_{pr},\T^+_{pr},\T^-_{pr})$ are abstract spectra.
 \end{theorem}

 \doc We first prove       that
 the bitopological space $(S'(X),\T^+,\T^-)$ is an abstract spectrum.
 For doing this it is enough to show that
 $S'(X)$ is a closed subset of the Cantor cube $\D^X$ (see \ref{th1}).

 Let $\{p_\s\in S'(X):\s\in \Sigma\}$ be a net in the Cantor cube $\D^X$
 converging to a point $p\in\D^X$. We have to prove that $p\in S'(X)$, i.e.
 that $p$ is an l-prime ideal in $(X,+,\times)$.

 Exactly as in the proof of \ref{th2}, we show that $a,b\in p$ implies that
 $a+b\in p$ and that if $a\times b\in p$ then   $a\in p$ or $b\in p$.

 Let $f_\s=e(p_\s)$ and
 $f=e(p)$ (see \ref{lm20} for the notation). Then the net
 $\{f_\s,\s\in\Sigma\}$ in $\D^X$ converges to $f$, i.e., for every $x\in X$,
 the net $\{f_\s(x),\s\in\Sigma\}$ in the discrete space $\D$ converges to
 $f(x)$.

 Let $a,b\in X$ and $a+b\in p$. Then $f(a+b)=1$.
 Hence there exists a $\s_0\in \Sigma$ such that
 $f_\s(a+b)=1$ for every $\s\ge\s_0$. Consequently, for every $\s\ge\s_0$, we have
 that $a+b\in p_\s$. Then, for every $\s\ge\s_0$,
 we get that $a\in p_\s$ and $b\in p_\s$,
 i.e. $f_\s(a)=1$ and $f_\s(b)=1$. This implies that $f(a)=1$ and $f(b)=1$,
 i.e. $a\in p$ and $b\in p$.

 Let $a,b\in X$ be such that $a\in p$ or $b\in p$. Suppose that
 $a\times b\not\in p$. Then $f(a\times b)=0$. Hence there exists
 a $\s_0\in\Sigma$ such that $f_\s(a\times b)=0$ for every $\s\ge \s_0$.
 This means that for every $\s\ge\s_0$, we have that $a\times b\not\in p_\s$.
 Consequently, $a\not\in p_\s$ and $b\not\in p_\s$ for every $\s\ge \s_0$.
 We obtain that $f_\s(a)=0$ and $f_\s(b)=0$ for every $\s\ge\s_0$.
 This implies that $f(a)=0$ and $f(b)=0$, i.e. $a\not\in p$ and
 $b\not\in p$, which is a contradiction. Therefore, $a\times b\in p$.
 Hence, $p$ is an l-prime ideal in $(X,+,\times)$.
 This shows that $S'(X)$ is a closed subset
 of the Cantor cube $\D^X$. Hence,
 the bitopological space $(S'(X),\T^+,\T^-)$ is an abstract spectrum.

 If the prime ideals $p_\s$ in the above proof were proper, then,
 obviously, $p$ would be  also  proper.
 This shows that the set $S'(X)_{pr}$ is also
 a closed subset of the Cantor cube $\D^X$. So, the bitopological space
 $(S'(X)_{pr},\T^+_{pr},\T^-_{pr})$ is an abstract spectrum.
 \sq

 \begin{exa}\label{ex20}
 \rm
 Let $(L,\vee,\we)$ be a distributive lattice with $0$ and $1$ and let us put
 $x + y= x\vee y$ and
 $x\times y= x\we y$, for every $x,y\in L$.
 Then, applying the construction
 from \ref{th4} to the set $L$ with the operations $+$ and $\times$
 and with fixed points $0$ and $1$, we get
 the topological space $(S'(L)_{pr},\T^+_{pr})$. We assert that it coincides
 with the classical spectrum $spec(L)$
 of the distributive lattice $(L,\vee,\we)$.
 \end{exa}

 \doc
 We first prove that the sets $spec(L)$ and $S'(L)_{pr}$ coincide.

 Let $p\in S'(L)_{pr}$. Then $0\in p$ and $1\not\in p$. If $a,b\in p$ then
 $a + b\in p$ and, hence, $a\vee b\in p$. Let $c\in L$,
 $a\in p$ and $c\le a$.
 Then $c\vee a = a$, i.e. $c + a\in p$. Thus
 $c\in p$. If
 $c,d\in L$ and $c\we d\in p$ then $c\times d\in p$. Therefore
 $c\in p$ or $d\in p$. So, $p\in spec(L)$.

 Let $p\in spec(L)$. If $a,b\in p$ then $a\vee b\in p$, i.e.
 $a + b\in p$. Further, if
 $x,y\in L$ and $x + y\in p$, then
 $x\vee y\in p$ and $x\le x\vee y$, $y\le x\vee y$. Hence $x\in p$ and $y\in p$.
 So, $x+y\in p$ iff $x\in p$ and $y\in p$. Now, let $a\in p$ or $b\in p$.
 Then $a\we b\le a$ and $a\we b\le b$. Therefore $a\we b\in p$, i.e.
 $a\times b\in p$. Finally, if $x,y\in L$ and $x\times y\in p$ then
 $x\we y\in p$ and, hence,
 $x\in p$ or $y\in p$. So, $x\times y\in p$ iff $x\in p$ or $y\in p$.
 Since $0\in p$ and $1\not\in p$, we obtain that $p\in S'(L)_{pr}$.
 Therefore, we proved that $S'(L)_{pr}=spec(L)$.

 The proof of the equality $\T^+_{pr}=\cal O$ is analogous to the proof of the
 corresponding statement about $S(L)_{pr}$, given in the proof of \ref{ex2}.
 \sq


 \subsection{The Main Theorem}

 The main theorem of Section 2 -- Theorem {\ref{th5} below --
 will be proved here. For doing this we need some preliminary definitions
 and results.

 \begin{defi}\label{de6}
 \rm
 Let $(S,\T^+,\T^-)$ be an abstract spectrum. For every two points
 $a,b\in S$ we put $a\le b$ iff $cl_{(S,\T^-)}\{a\}\sbe cl_{(S,\T^-)}\{b\}$
 (i.e., $a\le b$ iff $a$ is a specialization of $b$ in the  topological
 space $(S,\T^-)$).
 \end{defi}

 \begin{rem}\label{re3}
 \rm
 (a) The relation $\le$ defined in \ref{de6} is a partial order on $S$
 since $(S,\T^-)$ is a $T_0$-space (see \ref{pr1}) and,
 as it is well known, the specialization is a partial order
 on every $T_0$-space.

 \smallskip

 \noindent(b) It is obvious that $a\le b$ iff $a\in cl_{(S,\T^-)}\{b\}$ iff
 $b\in cl_{(S,\T^+)}\{a\}$ iff
 $cl_{(S,\T^+)}\{b\}\sbe cl_{(S,\T^+)}\{a\}$.

 \smallskip

 \noindent(c) It is easy to see that if $a\in S$ then
 $cl_{(S,\T^+)}\{a\}=\{b\in S: b\ge a\}$ and
 $cl_{(S,\T^-)}\{a\}=\{b\in S: b\le a\}$.

 \smallskip

 \noindent(d) If the elements of an abstract spectrum $S$
 are prime (or l-prime) ideals
 defined as in Section 2.2 (i.e. $S=S(X)$, where $X$ is
 a set with two binary operations), then $a\le b$ iff $a\sbe b$, for
 $a,b\in S$.
 \end{rem}

 \begin{lm}\label{lm1}
 Let $(S,\T^+,\T^-)$ be an abstract spectrum. If
 the net $\{a_{\s},\s\in\Sigma\}$ converges to $a$ in $(S,\T^-)$,
 the net $\{b_{\s},\s\in\Sigma\}$ converges to $b$ in $(S,\T^+)$ and
 $a_{\s}\le b_{\s}$ for every $\s\in\Sigma$, then $a\le b$.
 \end{lm}

 \doc Let $U\in\L^+$ and $b\in U$. Then there exists a
 $\s_0\in\Sigma$ such that $b_\s\in U$ for every $\s\ge\s_0$.
 Suppose that $a\not\in U$.
 Then $S\stm U\in\T^-$ and $a\in S\stm U$. Hence there exists a
 $\s_1\in\Sigma$ such that $a_\s\in S\stm U$ for every $\s\ge\s_1$.
 Putting $\s'=\sup\{\s_0,\s_1\}$, we obtain that $b_{\s'}\in U$ and
 $a_{\s'}\not\in U$. Therefore $b_{\s'}\not\in cl_{(S,T^+)}\{a_{\s'}\}$, i.e.
 $a_{\s'}\not\le b_{\s'}$ -- a contradiction. Hence $a\in U$.
 This shows that
 $b\in cl_{(S,\T^+)}\{a\}$, i.e. $a\le b$.
 \sq

 \begin{lm}\label{lm2}
 Let $(S,\T^+,\T^-)$ be an abstract spectrum. If
 $A\sbe S$ and $(A,\le)$ is a directed set (where $\le$ is
 the restriction to $A$ of the partial order defined in \ref{de6}),
 then the set $A$ has supremum in the ordered set $(S,\le)$.
 \end{lm}

 \doc Since $(A,\le)$ is a directed set and $A\sbe S$,
 $\{a, a\in A\}$ is a net in the compact Hausdorff
 space $(S,\T)$ (where $\T=\sup\{\T^+,\T^-\}$) (see \ref{pr6})
 and, hence, it has a cluster point $b\in S$.
 We shall prove that $b=\sup\{a:a\in A\}$ in $(S,\le)$.
 Indeed, let $U\in \T^+$ and $b\in U$. Then $U\in\T$ and for every
 $a\in A$ there exists an $a'\in A$ such that $a'\ge a$ and $a'\in U$.
 Hence $A\sbe U$. This shows that $b\in cl_{(S,\T^+)}\{a\}$ for every
 $a\in A$, i.e. $b\ge a$ for every $a\in A$. Let now $b'\in S$ and $b'\ge a$
 for every $a\in A$. The point $b$ is a limit in $(S,\T)$ (and, hence, in
 $(S,\T^-)$)
 of a net $\{a_\s,\s\in\Sigma\}$
 that is finer than the net
 $\{a,a\in A\}$. Put $b_\s=b'$ for every $\s\in\Sigma$.
 Then the net $\{b_\s,\s\in\Sigma\}$ converges to $b'$ in $(S,\T^+)$.
  Since $a_\s\le b_\s$ for every $\s\in\Sigma$, we obtain,
  using \ref{lm1}, that $b\le b'$. Hence, $b=\sup A$.
 \sq

 \begin{lm}\label{lm3}
 Let $(S,\T^+,\T^-)$ be an abstract spectrum. If
 $A\sbe S$ and $(A,\le')$ is a directed set, where $\le'$ is the inverse to
 the restriction to $A$ of the partial order defined in \ref{de6}
 (i.e. $a'\le'a''$ iff $a'\ge a''$, for $a',a''\in A$),
 then the set $A$ has infimum in the ordered set $(S,\le)$.
 \end{lm}

 \doc The proof of this lemma is completely analogous to that
 of Lemma \ref{lm2}.
 \sq

 \begin{lm}\label{lm4}
 Let $(S,\T^+,\T^-)$ be an abstract spectrum.
 Then for every $s\in S$ there exists an $m\in S$ (resp. $m'\in S$)
 such that $s\le m$ (resp. $m'\le s$) and
 $m$ is a maximal (resp. $m'$ is a minimal)
 element of the ordered set $(S,\le)$ (where $\le$ is
 from \ref{de6}).
 \end{lm}

 \doc It follows from the Zorn lemma and \ref{lm2} (resp. \ref{lm3}).
 \sq

 \begin{nota}\label{no2}
 \rm
 Let $(S,\T^+,\T^-)$ be an abstract spectrum.
 We put $Max(S)=\{m\in S: m$   is   a   maximal
   element   of  $(S,\le)\}$ and
 $Min(S)=\{m\in S: m$   is   a   minimal   element
   of  $(S,\le)\}$
 (where $\le$ is from \ref{de6}). We shall denote by
 $\T^+_M$
 (resp. $\T^-_M$) the induced by $\T^+$ (resp. $\T^-$) topology on $Max(S)$,
 and by $\T^+_m$
 (resp. $\T^-_m$) the induced by $\T^+$ (resp. $\T^-$) topology on $Min(S)$.
 \end{nota}

 \begin{pro}\label{pr7}
 Let $(S,\T^+,\T^-)$ be an abstract spectrum.
 Then:

 \smallskip

 \noindent(a)   $(Max(S),\T^+_M)$ and $(Min(S),\T^-_m)$ are compact $T_1$-spaces;

 \smallskip

 \noindent(b) $(Min(S),\T^+_m)$ and $(Max(S),\T^-_M)$
 are $T_2$-spaces;

 \smallskip

 \noindent(c) $Min(S)$ is
 dense in $(S,\T^+)$ and $Max(S)$ is dense in $(S,\T^-)$.
 \end{pro}

 \doc  (a) We first prove that $(Max(S),\T^+_M)$ is a compact $T_1$-space.
 Since, for every $a\in S$, $cl_{(S,\T^+)}\{a\}=\{b\in S:b\ge a\}$
 (see \ref{re3}(c)), we obtain that $(Max(S),\T^+_M)$ is a $T_1$-space.
 Let $\{a_\s,\s\in\Sigma\}$ be a net in $(Max(S),\T^+_M)$. Then
 $\{a_\s,\s\in\Sigma\}$ is a net in the compact space $(S,\T^+)$
 (see \ref{pr1}) and, hence, it has a cluster point $a\in S$ in $(S,\T^+)$.
 Now, we can find a net $\{a_{\s'},\s'\in\Sigma'\}$ in $(Max(S),\T^+_M)$
 which is finer than the net $\{a_\s,\s\in\Sigma\}$ and converges to $a$ in
 $(S,\T^+)$. By \ref{lm4}, there exists an $a'\in Max(S)$ such that $a\le a'$.
 Then $a'\in cl_{(S,\T^+)}\{a\}$ and, hence, the net
 $\{a_{\s'},\s'\in\Sigma'\}$ converges to $a'$ in $(Max(S),\T^+_M)$. This shows
 that the net $\{a_\s,\s\in\Sigma\}$ has a cluster point in $(Max(S),\T^+_M)$.
 Therefore, the space $(Max(S),\T^+_M)$ is compact.

 The proof of the fact that $(Min(S),\T^-_m)$ is a compact $T_1$-space is
 analogous.

 \smallskip

 \noindent(b) We first prove that $(Min(S),\T^+_m)$ is a Hausdorff space.
 Indeed, let $a,b\in Min(S)$ and $a\not= b$. Suppose that for any
 $U,V\in\L^+$ such that $a\in U$ and $b\in V$, we have that $U\cap V\not=\ems$.
 Then the family ${\cal F}=\{W\in\L^+:$
 $a\in W$ or $b\in W\}$ has the
 finite intersection property (see \ref{pr0}) and its elements are closed
 subsets of the compact space $(S,\T^-)$. Consequently there exists a
 $c\in\bigcap\cal F$. Since $\L^+$ is a base for $\T^+$, we obtain that
 $a\in cl_{(S,\T^+)}\{c\}$ and
 $b\in cl_{(S,\T^+)}\{c\}$. Hence
 $c\le a$ and $c\le b$. Having in mind that $a,b\in Min(S)$, we
 get that $c=a$ and $c=b$, i.e. $a=b$, which is  a contradiction. Therefore,
 $(Min(S),\T^+_m)$ is a Hausdorff space.

 Analogously, one proves that $(Max(S),\T^-_M)$ is a Hausdorff space.

 \smallskip

 \noindent(c) We first prove that
 $Min(S)$ is dense in $(S,\T^+)$.
 Indeed, let $x\in U\in\T^+$. By \ref{lm4}, there exists an $a\in Min(S)$
 such that $a\le x$. Then $x\in cl_{(S,\T^+)}\{a\}$. Hence $a\in U\cap Min(S)$.
 Therefore,
 $Min(S)$ is dense in $(S,\T^+)$.

 The proof of the fact that
 $Max(S)$ is dense in $(S,\T^-)$ is analogous.
 \sq

\medskip

 Let's recall the definitions of the coherent spaces and coherent maps:

 \begin{defi}\label{de7} {\rm (see, for example, \cite{J})}
 \rm
 Let $(X,\T)$ be a topological space.

 \smallskip

 \noindent(a) We shall denote by $KO(X,\T)$ (or, simply, by  $KO(X)$)
  the family of all
 compact open subsets of $X$.

 \smallskip

 \noindent(b) A closed subset $F$ of $X$ is
 called {\em irreducible} if the equality
 $F=F_1\cup F_2$, where $F_1$ and $F_2$ are closed subsets of $X$, implies
 that   $F=F_1$ or $F=F_2$.

 \smallskip

 \noindent(c) We say that the space $(X,\T)$ is {\em sober} if it is a $T_0$-space and
 for every non-void irreducible subset $F$ of $X$ there exists a $x\in X$ such
 that $F=cl_X\{x\}$.

 \smallskip

 \noindent(d) The space $(X,\T)$ is called {\em coherent} if it is a compact sober
 space and the family $KO(X,\T)$ is a
 closed under finite intersections
 base for the topology  $\T$.

 \smallskip

 \noindent(e) A continuous map $f:(X',\T')\lra (X'',\T'')$ is called {\em coherent}
 if $U''\in KO(X'')$ implies that $f^{-1}(U'')\in KO(X')$.
 \end{defi}

 \begin{nota}\label{no3}
 \rm
 We denote by \CohSp\ the category of all coherent spaces and
 all coherent maps between them.
 \end{nota}

 \begin{theorem}\label{th5}
 The categories \SS\ and \CohSp\
 are isomorphic.
 \end{theorem}

 \doc We shall construct two covariant functors
 $F:\SS\lra\CohSp$ and $G:\CohSp\lra\SS$ such that
 $F\circ G=Id_{\CohSp}$ and $G\circ F=Id_{\SS}$.

 For every $(S,\T^+,\T^-)\in |\SS |$, we put
 $F(S,\T^+,\T^-)=(S,\T^+)$.
 We shall prove that $(S,\T^+)\in |\CohSp|$. Indeed, we have: a) the space
 $(S,\T^+)$ is compact (by \ref{pr1}); b) $KO(S,\T^+)=\L^+$ (by \ref{pr2}) and
 hence the family $KO(S,\T^+)$ is a
 closed under finite intersections
 base for the topology  $\T^+$
 (by \ref{pr0} and
 (SP1) of \ref{de1}).
 Therefore we need only to show
 that $(S,\T^+)$ is a sober space. We have that
 $(S,\T^+)$ is a $T_0$-space (by \ref{pr1}). Let $A$ be a non-empty
 irreducible subset of $(S,\T^+)$. Then $A$ is a closed subset of $(S,\T)$,
 where $\T=\sup\{\T^+,\T^-\}$. Hence, by \ref{co2}, $(A,\T^+_A,\T^-_A)$
 is an abstract spectrum (where $\T^+_A$ (resp. $\T^-_A$) is the induced by
 $\T^+$ (resp. $\T^-$) topology on the subset
 $A$ of $S$). We shall prove
 that $|Min(A)|=1$. Suppose that $x,y\in Min(A)$ and $x\not=y$.
 Let $\T'$ be the induced by $\T^+_A$ topology on $Min(A)$. Since
 $(Min(A),\T')$ is a Hausdorff space (by \ref{pr7}(b)),
 there exists an $U\in\T'$
 such that $x\in U$ and $y\not\in cl_{(Min(A),\T\/')}U$.
 Put $B=cl_{(Min(A),\T\/')}U$ and $C=Min(A)\stm U$. Then $B$ and $C$ are
 closed subsets of $(Min(A),\T')$,
 $Min(A)=B\cup C$, $B\not=Min(A)$ and $C\not=Min(A)$. Since $Min(A)$ is
 dense in $(A,\T^+_A)$
 (by \ref{pr7}(c)), we obtain that $A=B'\cup C'$, where
 $B'=cl_{(A,\T^+_A)}B$ and $C'= cl_{(A,\T^+_A)}C$. The sets $B'$ and $C'$
 are closed in $(S,\T^+)$ since they are closed in
 $(A,\T^+_A)$ and $A$ is closed in $(S,\T^+)$. Moreover, $B'\not=A$ and
 $C'\not=A$, because $B'\cap Min(A)=B$ and $C'\cap Min(A)=C$.
 Since $A$ is irreducible, we get a contradiction. Therefore,
 $|Min(A)|=1$. Let $Min(A)=\{a\}$.
 Then \ref{pr7}(c) implies that $A=cl_{(S,\T^+)}\{a\}$. So, $(S,\T^+)$ is
 a sober space. We proved that $(S,\T^+)$ is a coherent space.

 Let $f\in\SS((S_1,\T^+_1,\T^-_1),(S_2,\T^+_2,\T^-_2))$.
 We denote by
 $F(f):S_1\lra S_2$ the function defined by
 $F(f)(x)=f(x)$ for every $x\in S_1$. We shall show that
 $F(f):(S_1,\T^+_1)\lra(S_2,\T^+_2)$
 is a coherent map. Indeed, since $f$ is a \SS-morphism, we have that
 $F(f):(S_1,\T^+_1)\lra(S_2,\T^+_2)$
 is a continuous map. Let $K\sbe S_2$, $K\in\T^+_2$ and $K$ be a compact
 subspace of $(S_2,\T^+_2)$. Then, by \ref{pr2}, $K\in\L^+_2$, i.e.
 $S_2\stm K\in\T^-_2$. Hence $f^{-1}(K)\in\T^+_1$ and
 $f^{-1}(S_2\stm K)\in\T^-_1$. Since $S_1\stm f^{-1}(K)=f^{-1}(S_2\stm K)$,
 we obtain that $f^{-1}(K)\in\L^+_1$. Consequently, by \ref{pr2}, $f^{-1}(K)$ is a
 compact subspace of $(S_1,\T^+_1)$. So, we  proved that
 $F(f)\in\CohSp(F(S_1,\T^+_1,\T^-_1),F(S_2,\T^+_2,\T^-_2))$.
 The definition of $F(f)$ implies immediately that $F$ preserves
 the identity maps and that $F(f\circ g)=F(f)\circ F(g)$. Therefore, we
 constructed a functor $F:\SS\lra\CohSp$.

 Let now $(S,\T^+)\in |\CohSp |$, ${\cal B}^+=KO(S,\T^+)$
 and ${\cal B}^-=\{S\stm U:U\in {\cal B}^+\}$.
 Since ${\cal B}^+$ is closed under finite intersections and finite unions,
 we obtain that ${\cal B}^-$ has the same properties.
 Obviously, $\bigcup{\cal B}^-=S$.
 Hence the family $\T^-$ of all subsets of $S$ that are unions of subfamilies
 of ${\cal B}^-$ is a topology on $S$ and ${\cal B}^-$ is a
 base for the topological space $(S,\T^-)$.
 We shall show that
 the bitopological space $(S,\T^+,\T^-)$ is an abstract spectrum
 and we will put $G(S,\T^+)=(S,\T^+,\T^-)$.

 It is easy to see that ${\cal B}^+\sbe\L^+$ and ${\cal B}^-\sbe\L^-$
 (see \ref{no1} for the notation).
 Since,
 by the definition of a coherent space,
 the family ${\cal B}^+$ is a base for the topological space $(S,\T^+)$
 and since the family ${\cal B}^-$ is a base for the space $(X,\T^-)$,
 we obtain that $\L^+$ (resp. $\L^-$)
 is a base for $(S,\T^+)$ (resp. $(S,\T^-)$). Hence the condition (SP1) of
 \ref{de1} is fulfilled. The condition (SP3) of \ref{de1} is also fulfilled
 since $(S,\T^+)$ is a $T_0$-space. Let us put $\T=\sup\{\T^+,\T^-\}$.
 We shall prove that the space $(S,\T)$ is compact. This will imply immediately
 that the condition (SP2) of \ref{de1} is fulfilled.

Obviously, for proving that $(S,\T)$ is compact, it is enough to show that every cover of $S$ of the type $\Omega=\OM^+\cup\OM^-$, where $\OM^+$ (resp. $\OM^-$) is
a subfamily of $\BB^+\stm\{S\}$ (resp. $\BB^-\stm\{S\}$), has a finite subcover. Let $\OM^*$ be the family of all finite unions of the elements of $\OM^-$. Then
$\OM^*\sbe\BB^-$, $\bigcup\OM^-= \bigcup\OM^*$ and $(\OM^*,\sbe)$ is a directed set (i.e. for every $U,V\in\OM^*$ there exists a $W\in\OM^*$ such that $U\cup V\sbe
W$). Put $H=S\stm\bigcup\OM^+$. Then $H\sbe\bigcup\OM^*$ and $H$ is a closed and, hence, compact subset of $(S,\T^+)$. If we find a $U_0\in\OM^*$ such that $H\sbe
U_0$ then we will have that $S\stm U_0\sbe S\stm H=\bigcup\OM^+$. From $U_0\in\BB^-$ we will get that $S\stm U_0\in\BB^+$ and, hence, $S\stm U_0$ will be a compact
subset of $(S,\T^+)$ covered by $\OM^+$. Consequently there will be a finite subfamily $\OM^+_f$ of $\OM^+$ covering $S\stm U_0$. Then $\OM^+_f\cup\{U_0\}$ will
cover $S$. Therefore, we will find a finite subcover of $\OM$. So, it is enough to prove that there exists an $U_0\in\OM^*$ such that $H\sbe U_0$.

Put $\HH^+=\{V\cap H:V\in\BB^+\}$. Then $\HH^+$ is a base for the subspace $H$ of $(S,\T^+)$, $\HH^+$ is closed under finite unions and finite intersections,
$\HH^+$ is a distributive lattice with respect to the operations $\cup$ and $\cap$ and, since $H$ is closed in $(S,\T^+)$, all elements of $\HH^+$ are compact
subsets of $(S,\T^+)$. Furthermore, for every $U\in\OM^*$ we put $U^+=S\stm U$. Then $U^+\in\BB^+$ for every $U\in\OM^*$.

Suppose that for every $U\in\OM^*$ we have that $H\stm U\not=\ems$. Then $H\cap U^+\not=\ems$ for every $U\in\OM^*$. Since for every $U,V\in\OM^*$ there exists a
$W\in\OM^*$ such that $W^+\sbe U^+\cap V^+$, the family $\{H\cap U^+:U\in\OM^*\}$ has the finite intersection property. Hence it generates a filter $\p$ in $\HH^+$.
Let $\Phi$ be an ultrafilter in $\HH^+$ containing $\p$ and let $L=\bigcap\{cl_{(S,\T^+)}W:W\in\Phi\}$. Then $L$ is a non-empty closed subset of $(S,\T^+)$ and
$L\sbe H$. Moreover, $L\cap W_0\not=\ems$ for every $W_0\in\Phi$. Indeed, let $W_0\in\Phi$. Then $W_0\in\HH^+$ and, hence, $W_0$ is a compact subset of $(S,\T^+)$.
It is easy to see that the family $\{cl_{W_0}(W_0\cap W):W\in\Phi\}$ has the finite intersection property. Consequently $\ems\not=\bigcap \{cl_{W_0}(W_0\cap
W):W\in\Phi\}= W_0\cap\bigcap\{cl_H(W_0\cap W):W\in\Phi\}\sbe W_0\cap\bigcap\{cl_HW:W\in\Phi\}=W_0\cap L$. So, we  proved that $L\cap W_0\not=\ems$ for every
$W_0\in\Phi$. We shall prove now that $L$ is an irreducible subset of $(S,\T^+)$. Indeed, suppose that $L=A\cup B$, where $A$ and $B$ are closed subsets of
$(S,\T^+)$ and $A\not= L$, $B\not=L$. Then $(H\stm A)\cap L\not=\ems$ and $(H\stm B)\cap L\not=\ems$. Let $x\in (H\stm A)\cap L$. Then there exists a $W'\in\HH^+$
such that $x\in W'\sbe H\stm A$. Since $x\in L$, we obtain that $W'\cap W\not=\ems$ for every $W\in\Phi$. Consequently $W'\in\Phi$. Analogously, taking an $y\in
(H\stm B)\cap L$, we can find a $W''\in\Phi$ such that $y\in W''\sbe H\stm B$. Putting $W_0=W'\cap W''$, we get that $W_0\in\Phi$. Since $W_0\sbe (H\stm A)\cap
(H\stm B)=H\stm(A\cup B)=H\stm L$, we conclude that $W_0\cap L=\ems$ -- a contradiction. Therefore, $L$ is an irreducible subset of $(S,\T^+)$. This implies,
because of the fact that $(S,\T^+)$ is sober, that there exists a point $l\in L$ such that $L=cl_{(S,\T^+)}\{l\}$. We shall show that
$l\in\bigcap\{U^+:U\in\OM^*\}$. Indeed, let $U\in\OM^*$. Then $H\cap U^+\in\p\sbe\Phi$. Hence $U^+\cap L\not=\ems$. Let $x\in U^+\cap L$. Then $x\in U^+\in\T^+$ and
$x\in L=cl_{(S,\T^+)}\{l\}$. Consequently $l\in U^+$. So, we  proved that $l\in\bigcap\{U^+:U\in\OM^*\}$. On the other hand we have that $l\in L\sbe
H\sbe\bigcup\OM^*= \bigcup\{S\stm U^+:U\in\OM^*\}=S\stm \bigcap\{U^+:U\in\OM^*\}$, i.e. $l\not\in\bigcap\{U^+:U\in\OM^*\}$ -- a contradiction. It shows that there
exists a $U_0\in\OM^*$ such that $H\sbe U_0$. Therefore, we  proved that the space $(S,\T)$ is compact and, hence,
 that the condition (SP2) of \ref{de1} is fulfilled. So,
 the bitopological space $(S,\T^+,\T^-)$ is an abstract spectrum.

 Let $f\in\CohSp((S_1,\T^+_1),(S_2,\T^+_2))$.
 We denote by
 $G(f):S_1\lra S_2$ the function defined by $G(f)(x)=f(x)$ for every
 $x\in S_1$. We shall show that
 $G(f)\in\SS((S_1,\T^+_1,\T^-_1),(S_2,\T^+_2,\T^-_2))$,
 where $(S_i,\T^+_i,\T^-_i)=G(S_i,\T^+_i)$, $i=1,2$.
 Indeed, we have that
 $f:(S_1,\T^+_1)\lra(S_2,\T^+_2)$
 is a continuous map and hence
 $G(f):(S_1,\T^+_1)\lra(S_2,\T^+_2)$
 is a continuous map. For proving that
 $G(f):(S_1,\T^-_1)\lra(S_2,\T^-_2)$
 is a continuous map it is enough to show that $U\in\BB^-_2$ implies that
 $f^{-1}(U)\in\BB^-_1$ (because $\BB^-_1$ (resp. $\BB^-_2$) is a base for
 $\T^-_1$ (resp. $\T^-_2$))
 (here we use the notation introduced above in the process of
 the definition of $G$ on the objects of the category \CohSp).
 So, let $U\in\BB^-_2$. Then
 $S_2\stm U\in KO(S_2,\T^+_2)$. Since $f$ is a coherent map, we obtain that
 $V=f^{-1}(S_2\stm U)\in KO(S_1,\T^+_1)=\BB^+_1$.
 Obviously, $V=S_1\stm f^{-1}(U)$.
 Consequently $f^{-1}(U)=S_1\stm V\in\BB^-_1$. So,
 $G(f)\in\SS(G(S_1,\T^+_1),G(S_2,\T^+_2))$.
 The definition of $G(f)$ implies immediately that $G$ preserves
 the identity maps and $G(f\circ g)=G(f)\circ G(g)$. Therefore, we
 constructed a functor $G:\CohSp\lra\SS$.

 From \ref{co1} and the constructions of the functors $F$ and $G$
 we get that
 $F\circ G=Id_{\CohSp}$ and $G\circ F=Id_{\SS}$.
 So, the categories \SS\ and \CohSp\ are isomorphic.
 \sq

 \begin{cor}\label{co3}
 The categories \DLat\ and \SS\ are dual.
 \end{cor}

 \doc Since the categories
 \DLat\ and \CohSp\ are dual (see, for example, \cite{J}),
 our statement follows immediately from \ref{th5}.
 \sq

 \begin{nist}\label{ni1}
 \rm
 Let us recall the descriptions of the
 duality functors
 $$F':\CohSp\lra\DLat\ \ \mbox{ and }\ \
 G':\DLat\lra\CohSp$$
  (see, for example, \cite{J}):
 if $(X,\T^+)$ is a coherent
 space then $$F'(X,\T^+)=(KO(X,\T^+),\cup,\cap,\ems,X);$$
  if
 $f\in\CohSp((X_1,\T^+_1),(X_2,\T^+_2))$ then $F'(f):F'(X_2,\T^+_2)\lra
 F'(X_1,\T^+_1)$ is defined by the formula $$F'(f)(U)=f^{-1}(U)$$ for every
 $U\in KO(X_2,\T^+_2)$; if $(L,\vee,\we,0,1)\in |\DLat |$ then
 $$G'(L,\vee,\we,0,1)=(spec(L),{\cal O}),$$ where $\cal O$ is the Stone
 topology on $spec(L)$ (see the proof of \ref{ex2} for
 the notation); if
 $f\in\DLat((L_1,\vee_1,\we_1,0_1,1_1),(L_2,\vee_2,\we_2,0_2,1_2))$
 then $$G'(f):
 G'(L_2,\vee_2,\we_2,0_2,1_2)\lra G'(L_1,\vee_1,\we_1,0_1,1_1))$$
 is defined by the formula $$G'(f)(p)=f^{-1}(p)$$ for every $p\in spec(L_2)$.
 The natural equivalence $\psi:Id_{\CohSp}\lra G'\circ F'$ is given by the
 formula $\psi(X,\T^+)=\psi_{(X,\T^+)}$
 for every $(X,\T^+)\in |\CohSp |$, where $$\psi_{(X,\T^+)}:(X,\T^+)
 \lra (G'\circ F')(X,\T^+),\ \ x\mapsto \{U\in F'(X,\T^+):x\not\in U\}.$$
 In particular, $\psi_{(X,\T^+)}$ is a \CohSp -isomorphism for every
 coherent space $(X,\T^+)$. The natural equivalence
 $\phi:Id_{\DLat}\lra F'\circ G'$ is given by the formula
 $\phi(L)=\phi_L$ for every $L\in |\DLat |$, where
 $$\phi_L:L\lra (F'\circ G')(L),\ \  l\mapsto \{p\in G'(L):l\not\in p\}.$$
 In particular, $\phi_L$ is a \DLat -isomorphism for every distributive
 lattice $L$.
 \end{nist}


 \subsection{Some Applications}

 Let's start with recalling that if $L$ is a distributive lattice with
 $0$ and $1$ then its classical spectrum $spec(L)$ can be interpreted
 as an abstract spectrum (see \ref{ex2}, \ref{pr3} and \ref{co1}).

 We will first prove a general theorem.

 \begin{theorem}\label{th66}
 Let $X$ be a set, $S$ be a family of subsets of $X$ (i.e. $S\sbe Exp(X)$),
 $\T^+$ and $\T^-$ be the topologies on $S$ defined in \ref{th1}, and let the
 bitopological space $(S,\T^+,\T^-)$ be an abstract spectrum.
 Then there exist a distributive lattice $L$ with 0 and 1, and a function
 $\p : X\lra L$ such that:

 \smallskip

 \noindent(i) the set $\p(X)$ generates $L$;

 \smallskip

 \noindent(ii) $\p^{-1}(q)\in S$ for every $q\in spec(L)$
 (see \ref{ex2} for the notation);

 \smallskip

 \noindent(iii) $\Phi: spec(L)\lra S$, $q\mapsto \p^{-1}(q)$,
 is an \SS-isomorphism;

 \smallskip

 \noindent(iv) if $L'$ is a distributive lattice with 0 and 1, and
 $\t :X\lra L'$ is a function such that:

 \hspace{0.5cm}$(1)$ $\t^{-1}(q)\in S$
 for every $q\in spec(L')$, and

 \hspace{0.5cm}$(2)$ $\Theta:spec(L')\lra S$, $q\mapsto \t^{-1}(q)$,
 is an \SS-morphism,\\
 then there exists a unique lattice homomorphism
 $l:L\lra L'$ with $l\circ\p=\t$;

 \smallskip

 \noindent(v) if $\p_1:X\lra L_1$,
 where $L_1$ is a distributive lattice with 0 and 1, is
 such that:

 \hspace{0.5cm}$(1')$ $(\p_1)^{-1}(q)\in S$ for every $q\in spec(L_1)$, and

 \hspace{0.5cm}$(2')$ $\Phi_1: spec(L_1)\lra S$, $q\mapsto (\p_1)^{-1}(q)$,
 is an \SS-isomorphism,\\
 then there exists a unique lattice isomorphism $l:L\lra L_1$ with
 $l\circ \p = \p_1$;

 \smallskip

 \noindent(vi) $\p:X\lra L$ is an injection iff for any two different
 points $x$ and $y$ of $X$ there exists a $p\in S$ containing exactly one
 of them.
 \end{theorem}

 \doc We shall use the notation of \ref{th1}, \ref{th2} and \ref{ex2}.

 By  (the proof of) \ref{th5}, we have that
 $(S,\T^+)\in |\CohSp |$. We
 put $L=F'(S,\T^+)$ (see \ref{ni1}), i.e.
 $L=\{U\in\T^+:U$ is compact$\}$ and, hence, by \ref{pr2}, $L=\L^+$.
 Then $L$ is a distributive lattice
 with 0 and 1. Define the function $\p:X\lra L$ by the formula
 $\p(x)=U^+_x$ for every $x\in X$ (recall that
 $U^+_x=\{p\in S:x\not\in p\}$ and $U^+_x\in\L^+$ (see \ref{th1}
 and the part $(b)\Ra (a)$ of its proof)). Hence
  $\p(X)$
 ($=\{U^+_x:x\in X\} =
 \PP^+$)
 is a subbase for $\T^+$
 (see \ref{th1}).
In what follows,
 the topological space $(S,\T^+)$ will be denoted, briefly, by
 $S$.

 \smallskip

 \smallskip

 \noindent {\em The proof of (i)}:
 Let $L^*$ be the set of all finite unions
 of the elements of the set $\BB^+$
 of all finite intersections of the elements of
 $\PP^+=\p(X)$. Then $L^*$ coincides with the subset of $L$ generated by
 $\p(X)$ and $\BB^+$ is a  base for $\T^+$.
 If $U\in L$ then $U$ is a
 compact open subset of $S$ and, hence, it is a finite union of
 elements of $\BB^+$. Thus $U\in L^*$. Therefore, the set $\p(X)$ generates
 $L$.

 \smallskip

 \smallskip

 \noindent{\em The proof of (ii) and (iii)}:
 By \ref{ni1}, we have that $spec(L)=G'(L)$.
 Since the map $\psi_S:S\lra (G'\circ F')(S)$,
 $p\lra \{U\in L:p\not\in U\}$ is a \CohSp -isomorphism (see \ref{ni1}),
 we get that $spec(L)=\psi_S(S)$.

 Let $q\in spec(L)$. Then there exists a unique $p\in S$ such that
 $q=\psi_S(p)$.
 So, we have that
 $\p^{-1}(q)=\p^{-1}(\psi_S(p))=
 \{x\in X:\p(x)\in\psi_S(p)\}=
 \{x\in X:U^+_x\in\psi_S(p)\}=\{x\in X:p\not\in U^+_x\}=
 \{x\in X:x\in p\}=p$, i.e. $\p^{-1}(q)=\psi_S^{-1}(q)$ for every
 $q\in spec(L)$. Since the function $\psi_S^{-1}$ is a \CohSp -isomorphism,
 we conclude that the function
 $\Phi: spec(L)\lra S$, $q\lra \p^{-1}(q)$, is a \CohSp -isomorphism.
 Now, (the proof of) \ref{th5} implies, that $\Phi$ is an \SS -isomorphism.

 \smallskip

 \smallskip

 \noindent{\em The proof of (iv)}:
 Put $\tau=\psi_S\circ \TE$. Then, by \ref{th5} and \ref{ni1},
 $\TE: spec(L')\lra (S,\T^+)$ and $\tau: spec(L')\lra (G'\circ F')(S,\T^+)$
 are \CohSp -morphisms. Since $G'(L')=spec(L')$  and
 $F'(S,\T^+)=L$, we obtain that
 $F'(\tau)=F'(\TE)\circ F'(\psi_S):(F'\circ G')(L)\lra (F'\circ G')(L')$
 (see \ref{ni1}). Put
 $l=\phi_{L'}^{-1}\circ F'(\tau)\circ \phi_L$
 (using the notation from \ref{ni1}).
 Then $l:L\lra L'$ is a
 lattice homomorphism. We shall prove that
 $F'(\TE)\circ F'(\psi_S)\circ \phi_L\circ\p=\phi_{L'}\circ\t$. This will
 imply that
 $\phi_{L'}^{-1}\circ F'(\TE)\circ F'(\psi_S)\circ \phi_L\circ\p=\t$ and,
 hence, we wll have that
 $\t=\phi_{L'}^{-1}\circ (F'(\TE)\circ F'(\psi_S))\circ \phi_L\circ\p=
 (\phi_{L'}^{-1}\circ  F'(\tau)\circ \phi_L)\circ\p= l\circ\p$, i.e. that
 $\t=l\circ\p$.

 Let $x\in X$. Then $(\phi_{L'}\circ\t)(x)=\phi_{L'}(\t(x))=\{q'\in spec(L'):
 \t(x)\not\in q'\}$. On the other hand, $(\phi_L\circ\p)(x)=\phi_L(\p(x))=
 \{q\in spec(L):\p(x)\not\in q\}$.
 Put $U=(F'(\psi_S)\circ \phi_L\circ\p)(x)$.
 Since $\psi_S^{-1}=\Phi$
 (see the proof of (ii) and (iii) above), we get that
 $(F'(\psi_S))^{-1}=F'(\psi_S^{-1})=F'(\Phi)$. Hence
 $(F'(\Phi))(U)=(F'(\psi_S))^{-1}(U)=(\phi_L\circ\p)(x)$. Now, the
 definition of  $F'(\Phi)$ (see \ref{ni1}) implies that
 $(F'(\Phi))(U)=\Phi^{-1}(U)$.
 Hence $\Phi^{-1}(U)=(\phi_L\circ\p)(x)$.
 Since $\Phi$ is an isomorphism (see (iii)), we get that
 $U=\Phi((\phi_L\circ\p)(x))=\Phi(\{q\in spec(L):\p(x)\not\in q\})=
 \{\Phi(q):q\in spec(L),\p(x)\not\in q\}=
 \{\p^{-1}(q):q\in spec(L),\p(x)\not\in q\}=
 \{\p^{-1}(q):q\in spec(L),x\not\in \p^{-1}(q)\}=
 \{p\in S:x\not\in p\}=U^+_x$, i.e $U=U^+_x$.
 Therefore,
 $(F'(\psi_S)\circ \phi_L\circ\p)(x)=U^+_x$.
 Then
 $(F'(\TE)\circ F'(\psi_S)\circ \phi_L\circ\p)(x)=
 (F'(\TE))((F'(\psi_S)\circ \phi_L\circ\p)(x))=
 (F'(\TE))(U^+_x)=\TE^{-1}(U^+_x)=
 \{q'\in spec(L'):\TE(q')\in U^+_x\}=
 \{q'\in spec(L'): \t^{-1}(q')\in U^+_x\}=
 \{q'\in spec(L'): x\not\in\t^{-1}(q')\}=
 \{q'\in spec(L'):\t(x)\not\in q'\}=(\phi_{L'}\circ\t)(x)$.
 So, we  proved that
 $\t=l\circ\p$.
 This, combined with the fact that $\p(X)$ generates $L$ (see (i)), proves
 the uniqueness of $l$.

 \smallskip

 \noindent{\it The proof of (v)}: Let $\p_1:X\lra L_1$
 has the properties $(1')$ and $(2')$.
 Then, using (iv), we obtain a lattice
 homomorphism $l:L\lra L_1$ such that $l\circ\p=\p_1$. From the construction
 of $l$, given in (iv), we have that
 $l=\phi_{L_1}^{-1}\circ F'(\psi_S\circ \Phi_1)\circ \phi_L$.
 Since $\Phi_1$ is an \CohSp -isomorphism (by $(2')$ and \ref{th5}), we
 get that $l$ is a \DLat -isomorphism (because all other components of the
 composition defining $l$ are also \DLat -isomor\-phisms (see \ref{ni1})).

 \smallskip

 \smallskip

 \noindent{\it The proof of (vi)}: Let $x,y\in X$
 and $x\not=y$. Then $\p(x)=\{p\in S:x\not\in p\}$
 and $\p(y)=\{p\in S:y\not\in p\}$. Hence, $\p(x)\not=\p(y)$ iff there exists
 a $p\in S$ containing exactly one of the points $x$ and $y$.
 \sq

 \begin{cor}\label{th6}
 Let $X$ be a set endowed with two arbitrary multivalued binary operations
 $\op$ and $\ot$ and
 with two fixed different points $\xi_0\in X$ and $\xi_1\in X$.
 Then there exist a distributive lattice $(L,\vee,\we)$
 with 0 and 1, and a function
 $\p : X\lra L$ such that:

 \smallskip

 \noindent(i) the set $\p(X)$ generates $L$;

 \smallskip

 \noindent(ii) $\p^{-1}(q)\in S(X)_{pr}$ for every $q\in spec(L)$
 (resp. $\p^{-1}(q)\in S(X)$ for every $q\in spec(L)$)
 (see \ref{th2} and \ref{ex2} for the notation);

 \smallskip

 \noindent(iii) $\Phi: spec(L)\lra S(X)_{pr}$, $q\mapsto \p^{-1}(q)$
 (resp. $\Phi: spec(L)\lra S(X)$, $q\lra \p^{-1}(q)$)
 is an \SS-isomorphism;

 \smallskip

 \noindent(iv) if $L'$ is a distributive lattice with 0 and 1, and
 $\t :X\lra L'$ is a function such that:

 \hspace{0.5cm}$(1)$ $\t^{-1}(q)\in S(X)_{pr}$
 (resp. $\t^{-1}(q)\in S(X)$)
 for every $q\in spec(L')$, and

 \hspace{0.5cm}$(2)$ $\Theta:spec(L')\lra S(X)_{pr}$, $q\mapsto \t^{-1}(q)$,
 (resp. $\Theta:spec(L')\lra S(X)$, $q\mapsto \t^{-1}(q)$,)
 is an \SS-morphism,\\
 then there exists a unique lattice homomorphism
 $l:L\lra L'$ with $l\circ\p=\t$;

 \smallskip

 \noindent(v) if $\p_1:X\lra L_1$, where $L_1$ is a distributive lattice with 0 and 1, is
 such that:

 \hspace{0.5cm}$(1')$ $(\p_1)^{-1}(q)\in S(X)_{pr}$ for every $q\in spec(L_1)$
 (resp. $(\p_1)^{-1}(q)\in S(X)$ for every $q\in spec(L_1)$), and

 \hspace{0.5cm}$(2')$ $\Phi_1: spec(L_1)\lra S(X)_{pr}$, $q\mapsto (\p_1)^{-1}(q)$
 (resp. $\Phi_1: spec(L_1)\lra S(X)$, $q\mapsto (\p_1)^{-1}(q)$)
 is an \SS-isomorphism,\\
 then there exists a unique lattice isomorphism $l:L\lra L_1$ with
 $l\circ \p = \p_1$;

 \smallskip

 \noindent(vi) $a\op b\sbe\{x\in X: \p(x)\le \p(a)\vee\p(b)\}$ and
 $a\ot b\sbe\{x\in X: \p(x)\ge \p(a)\we\p(b)\}$ for any $a,b\in X$.
 \end{cor}

 \doc  Denote by $S$ the set $S(X)_{pr}$ (resp. $S(X)$)
 (see \ref{th2} for the notation) and
 define the topologies $\T^+_{pr}$ (resp. $\T^+$) and $\T^-_{pr}$ (resp. $\T^-$)
 on $S$ as in \ref{th1}.
 Then, by \ref{th2}, the bitopological space $(S,\T^+_{pr},\T^-_{pr})$
 (resp. $(S,\T^+,\T^-)$)
 is an abstract
 spectrum. Hence, applying Theorem \ref{th66}, we obtain a distributive
 lattice $(L,\vee,\we,0,1)$ and a function
 $\p:X \lra L$ satisfying  conditions (i)-(v) of \ref{th66}
 and, hence, our conditions (i)-(v) as well.
 Consequently, we need
 only to check that the condition (vi) is also satisfied.
 In what follows, the notation of the proof of \ref{th66} and the
 construction of the function $\p$ given there are used.

 Let $a,b\in X$ and $x\in a\op b$. Then $\p(a)\vee\p(b)=\p(a)\cup\p(b)=
 \{p\in S:$   $a\not\in p$ or $b\not\in p\}$.
 Hence $S\stm (\p(a)\cup\p(b))=
 \{p\in S:a\in p$ and $b\in p\}$. Let $p'\in\p(x)=U^+_x=\{p\in S:x\not\in p\}$
 and suppose that
 $p'\not\in \p(a)\cup\p(b)$. Then $a\in p'$ and $b\in p'$. This implies that
 $a\op b\sbe p'$. Then $x\in p'$ and, hence, $p'\not\in\p(x)$ --
 a contradiction. Therefore, $p'\in\p(a)\cup\p(b)$. This shows that
 $\p(x)\sbe\p(a)\cup\p(b)$, i.e. $\p(x)\le\p(a)\vee\p(b)$, for every
 $x\in a\op b$. Consequently,
 $a\op b\sbe\{x\in X: \p(x)\le \p(a)\vee\p(b)\}$
 for any $a,b\in X$.

 Let $x\in a\ot b$.
 We have that $\p(a)\we\p(b)=\p(a)\cap\p(b)=\{p\in S:a\not\in p$ and
 $b\not\in p\}$. Let $p'\in\p(a)\cap\p(b)$.
 Then $a\not\in p'$ and $b\not\in p'$.
 Suppose that $p'\not\in\p(x)$. Then $x\in p'$ and, hence,
 $(a\ot b)\cap p'\not=\ems$.
 This implies that  $a\in p'$ or $b\in p'$, i.e. we get a  contradiction.
 Therefore, $p'\in\p(x)$. So, $\p(a)\cap\p(b)\sbe\p(x)$, i.e.
 $\p(a)\we\p(b)\le\p(x)$ for every $x\in a\ot b$.
 \sq

 \begin{cor}\label{th7}
 Let $X$ be a set endowed with two arbitrary single-valued binary operations
 $+$ and $\times$
 and with two fixed different points $\xi_o\in X$ and $\xi_1\in X$.
 Then there exist a distributive lattice
 $(L,\vee,\we)$ with 0 and 1, and a function
 $\p : X\lra L$ such that:

 \smallskip

 \noindent(i) the set $\p(X)$ generates $L$;

 \smallskip

 \noindent(ii) $\p^{-1}(q)\in S'(X)$ for every $q\in spec(L)$
 (resp. $\p^{-1}(q)\in S'(X)_{pr}$ for every $q\in spec(L)$)
 (see \ref{th4}, \ref{ex2} and \ref{th2} for the notation);

 \smallskip

 \noindent(iii) $\Phi: spec(L)\lra S'(X)$, $q\mapsto \p^{-1}(q)$,
 (resp. $\Phi: spec(L)\lra S'(X)_{pr}$, $q\mapsto \p^{-1}(q)$,)
 is an \SS-isomorphism;

 \smallskip

 \noindent(iv) if $L'$ is a distributive lattice with 0 and 1, and
 $\t :X\lra L'$ is a function such that:

 \hspace{0.5cm}$(1)$ $\t^{-1}(q)\in S'(X)$
 (resp. $\t^{-1}(q)\in S'(X)_{pr}$)
 for every $q\in spec(L')$, and

 \hspace{0.5cm}$(2)$ $\Theta:spec(L')\lra S'(X)$, $q\mapsto \t^{-1}(q)$
 (resp. $\Theta:spec(L')\lra S'(X)_{pr}$, $q\mapsto \t^{-1}(q)$)
 is an \SS-morphism,\\
 then there exists a unique lattice homomorphism
 $l:L\lra L'$ with $l\circ\p=\t$;

 \smallskip

 \noindent(v) if $\p_1:X\lra L_1$, where $L_1$ is a distributive lattice with 0 and 1, is
 such that:

 \hspace{0.5cm}$(1')$ $(\p_1)^{-1}(q)\in S'(X)$ for every $q\in spec(L_1)$
 (resp. $(\p_1)^{-1}(q)\in S'(X)_{pr}$ for every $q\in spec(L_1)$), and

 \hspace{0.5cm}$(2')$ $\Phi_1: spec(L_1)\lra S'(X)$, $q\mapsto (\p_1)^{-1}(q)$
 (resp. $\Phi_1: spec(L_1)\lra S'(X)_{pr}$, $q\mapsto (\p_1)^{-1}(q)$)
 is an \SS-isomorphism,\\
 then there exists a unique lattice isomorphism $l:L\lra L_1$ with
 $l\circ \p = \p_1$;

 \smallskip

 \noindent(vi) $\p(a+b) = \p(a)\vee\p(b)$ and $\p(a\times b)=\p(a)\we\p(b)$ for
 every $a,b\in X$.
 \end{cor}

 \doc  Denote by $S$ the set $S'(X)$ (resp. $S'(X)_{pr}$)
 (see \ref{th4} for the notation) and
 introduce the topologies $\T^+$ (resp. $\T^+_{pr}$) and $\T^-$
 (resp. $\T^-_{pr}$) on $S$ as in \ref{th1}.
 Then, by \ref{th4}, the bitopological space
 $(S,\T^+,\T^-)$
 (resp. $(S,\T^+_{pr},\T^-_{pr})$ )
 is an abstract
 spectrum. Hence, applying Theorem \ref{th66}, we obtain a distributive
 lattice $(L,\vee,\we,0,1)$ and a function
 $\p:X \lra L$ satisfying  conditions (i)-(v) of \ref{th66}
 and, hence, our conditions (i)-(v) as well. Consequently, we need
 only to check that the condition (vi) is also satisfied. This can be done
 easily (see the proof of \ref{th6}).
 \sq


\section{Separative algebras}

The main aim of this section is to give a detailed exposition of the theory of separative algebras, introduced and announces by Prodanov in \cite{P1}. This theory is a straight generalization of the theory of convex spaces in the sense of Tagamlitzki \cite{T}, which have been also a subject of Prodanov's Ph.D. dissertation \cite{PD}. We will follow very closely the style of  Prodanov's proofs from \cite{P5} and \cite{PD}.

\subsection{Preseparative algebras}

Let $X\not=\emptyset$ be a set with two binary multivalued operations denoted by $``\times"$ and $``+"$. This means that for any $x,y\in X$, $x \times y\subseteq X$ and $x + y\subseteq X$. Later on, instead of $``\times"$ and $``+"$, we shall use $``."$ and $``+"$, and following the common mathematical practice, sometimes we shall omit the sign $``."$.

We extend the operations $``."$ and $``+"$ for arbitrary subsets $A$ and $B$ of $X$ putting
$$A.B=\bigcup_{a\in A,b\in B} a.b\mbox{ and }A+B=\bigcup_{a\in A,b\in B} a+b$$
The one element subset $\{x\}\subseteq X$ will be denoted simply by $x$. Then for instance $x(yz)$ will mean $\{x\}.(y.z)$.

\begin{defi}\label{3.1.1}
\rm The system $\unlx=(X, .,+)$ is called a {\em preseparative
algebra} if $X\not=\emptyset$, $``."$ and $``+"$ are binary multivalued operations in X
satisfying the following axioms: for arbitrary $a,b,c,x\in X$, \\
(i)    $ab=ba$, \ \ \ \ \ \ \ \ \ \ \ \ \ \ \ \	(i')   $a+b=b+a$, \\
(ii)   $a(bc)=(ab)c$,   \ \ \ \ \ \ \ \ \  (ii')  $a+(b+c)=(a+b)+c$, \\
(iii)  from $a\in b+x$, and $c\in dx$, it follows that $(ad)\cap(b+c)\not=\emptyset$.
\end{defi}

By means of the operations  $``."$  and $``+"$, we introduce two new operations as follows: \\
division:\ \ \ \  $a/b=\{x\in X : a\in b.x\}$ and \\
difference: $a-b=\{x\in X : a\in b+x\}$.

We extend the operations division and difference for arbitrary subsets putting
$$A/B=\bigcup_{a\in A, b\in B}  a/b,\ \ \ \   A-B=\bigcup_{a\in A, b\in B}   a-b.$$
Sometimes instead of $A/B$ we will write $A:B$ or $\frac{A}{B}$.

The following lemma follows  immediately from the relevant definitions.

\begin{lm}\label{3.1.2} Let $``\bullet"$ be any of the operations $``."$, $``+"$, $``/"$ and $``-"$. Then the following conditions are true: \\
(i) $A\bullet\emptyset=\emptyset\bullet A=\emptyset$, \\
(ii) If $A\subseteq A'$ and $B\subseteq B'$ then $A\bullet B\subseteq A'\bullet B'$, \\
(iii) $(\bigcup_{i\in I} A_i)\bullet(\bigcup_{j\in J} B_j)=\bigcup_{i\in I, j\in J} A_i\bullet B_j$   and, in particular,  \\
(iii') $A\bullet (B\cup C)=(A\bullet B)\cup(A\bullet C)$, \\
(iv)  $(\bigcap_{i\in I} A_i)\bullet(\bigcap_{j\in J} B_j)\subseteq \bigcap_{i\in I,j\in J} A_i\bullet B_j$.
\end{lm}

\begin{pro}\label{3.1.3}  The following is true for arbitrary $A,B,C\subseteq X$: \\
(i)  $(A/B)\cap C\not=\emptyset$ iff $A\cap(B.C)\not=\emptyset$, \\
(ii) $(A-B)\cap C\not=\emptyset$ iff $A\cap(B+C)\not=\emptyset$.
\end{pro}

\doc
(i) $(A/B)\cap C\not=\emptyset$ iff $\exists x\in X$: $x\in (A/B)\cap C$ iff $\exists x\in X$: $x\in A/B$ and $x\in C$ iff $\exists x,a,b\in X$ $a\in A$, $b\in B$ $x\in a/b$ and $x\in C$ iff $\exists x,a,b\in X$: $a\in A$, $b\in B$, $a\in b.x$ and $x\in C$ iff $\exists a\in X$: $a\in A$ and $a\in B.C$ iff $\exists a\in X$: $a\in A\cap (B.C)$ iff $A\cap(B.C)\not=\emptyset$.

The proof of (ii) is similar. \sq

\begin{pro}\label{3.1.4} The following conditions are true for arbitrary subsets $A,B,C$ of $X$: \\
(i)   $AB=BA$,\ \ \ \ \ \ \ \ \ \ \ \ \ \ (i')  $A+B=B+A$, \\
(ii)  $A(BC)=(AB)C$,\ \ \  (ii') $A+(B+C)=(A+B)+C$.
\end{pro}

\doc As an example we shall verify (i) . The proof of the remaining conditions is similar.

$x\in AB$ iff $\exists a\in A$ $\exists b\in B$: $x\in ab$ iff (by commutativity of $``."$) $\exists a\in A$ $\exists b\in B$: $x\in ba$ iff $x\in BA$. \sq

Associativity enables us to write $A_1.A_2.\dots A_n$  and $A_1+A_2+\dots+A_n$ without parentheses.

We denote $A^n =A.A....A$ (n-times) and $nA=A+A+...+A$ (n-times), putting $A^1=1A=A$.

\begin{lm}\label{3.1.5}
(i)   $A^iA^j=A^{i+j}$, \\
(i') $iA+jA=(i+j)A$ \\
(ii)   $(A\cup B)^2=A^2\cup AB\cup B^2$, \\
$\ \ \ \ \ (A\cup B\cup C)^2=A^2\cup AB\cup AC\cup BC\cup C^2$ \\
(ii') $2(A\cup B)=2A\cup(A+B)\cup2B$. \\
$\ \ \ \ \ \ 2(A\cup B\cup C)=2A\cup(A+B)\cup(A+C)\cup(B+C)\cup2C$.
\end{lm}

\doc  (i) and (i') follow immediately from the definition, and
(ii) and (ii') follow from Lemma \ref{3.1.2}(iii') and commutativity. \sq

\begin{pro}\label{3.1.6} The following conditions are equivalent to the Axiom (iii) from the definition of  preseparative algebras (see Definition \ref{3.1.1}): \\
(i) $a+\frac{b}{c}\subseteq\frac{a+b}{c}$, \\
(ii) $a(b-c)\subseteq ab-c$.
\end{pro}

\doc As an example we shall show the equivalence of the Axiom (iii)  with (i).

\smallskip

\noindent((Axiom (iii))$\longrightarrow$ (i)). Let $x\in a+\frac{b}{c}$. Then there exists $y\in X$ such that $x\in a+y$, $y\in\frac{b}{c}$ and $b\in c+y$. By Axiom (iii), $(xc)\cap(a+b)\not=\emptyset$. Then, by Proposition \ref{3.1.3}(i), we obtain that $x\cap\frac{a+b}{c}\not=\emptyset$ and hence $x\in\frac{a+b}{c}$. Since $x$ is an arbitrary element of $X$, this shows that $a+\frac{b}{c}\subseteq\frac{a+b}{c}$.

\smallskip

\noindent((i)$\longrightarrow$  (Axiom  (iii))).  Let  $a\in b+x$ and $c\in dx$.  Then $x\in\frac{c}{d}$ and $c\in b+\frac{c}{d}$. Then, by (i),  $c\in\frac{b+c}{d}$, so
that $c\cap\frac{b+c}{d}\not=\emptyset$. Applying Proposition \ref{3.1.3}(i), we obtain that $(cd)\cap(b+c)\not=\emptyset$, which shows that Axiom (iii) holds.

The equivalence of Axiom (iii) with (ii) can be proved similarly by using Proposi\-tion \ref{3.1.3}(ii). \sq

\begin{pro}\label{3.1.7} The following conditions are true for arbitrary subsets $A,B,C,D$ of $X$: \\
(i) $A+\frac{B}{C}\subseteq\frac{A+B}{C}$, \\
(ii) $A(B-C)\subseteq AB-C$, \\
(iii)  $(A/B)/C=A/(BC)$, \\
(iv)  $(A-B)-C=A-(B+C)$, \\
(v)  $\frac{A}{B}+\frac{C}{D}\subseteq\frac{A+C}{B.D}$, \\
(vi) $(A-B)(C-D)\subseteq AC-(B+D)$.
\end{pro}

\doc  (i) and (ii) are extensions of (i) and (ii) of Proposition \ref{3.1.6} for arbitrary sets and follow directly from Proposition \ref{3.1.6}.

\smallskip

\noindent(iii) Let $x$ be an arbitrary element of $X$. Then, applying Proposition \ref{3.1.3}(i), we obtain that $x\in(A/B)/C$ iff $(A/B)/C\cap x\not=\emptyset$ iff $(A/B)\cap Cx\not=\emptyset$ iff $A\cap (BCx)\not=\emptyset$ iff $A\cap(BC)x\not=\emptyset$ iff $(A/(BC))\cap x\not=\emptyset$ iff $x\in A/(BC)$. Hence,  $(A/B)/C=A/(BC)$.

\smallskip

\noindent(iv) The proof can be done similarly by applying Proposition \ref{3.1.3}(ii).

\smallskip

\noindent(v) $\frac{A}{B}+\frac{C}{D}\subseteq\frac{A/B+C}{D}\subseteq\frac{(A+C)/B}{D}=\frac{A+C}{B.D}$. We have applied
two times (i) and then (iii).

\smallskip

\noindent(vi) The proof can be done similarly by applying two times (ii) and then (iv). \sq

\subsection{Filters and ideals in preseparative algebras}

\begin{defi}\label{3.2.1}
\rm Let $\underline{X}=(X,.,+)$ be a preseparative algebra. A subset $F\subseteq X$ is called a {\em filter} in $X$ if $F.F\subseteq F$. A subset $I\subseteq X$ is called an {\em ideal} in $X$ if $I+I\subseteq I$. A subset $F\subseteq X$ is called a {\em prime filter} in $X$ if $F$ is a filter and the complement $X\setminus F$ of $F$ is an ideal in $X$. Dually, a subset $I\subseteq X$ is called a {\em prime ideal} in $X$ if $I$ is an ideal and $X\setminus I$ is a filter in $X$.
\end{defi}

Obviously the empty set $\emptyset$ and the whole set $X$ are examples of a filter, ideal, prime filter and prime ideal. They are in some sense trivial examples. Nontrivial examples of filters and ideals will be given by the constructions $\mu(A)$ and $\alpha(A)$ below. Constructions of prime filters and prime ideals will be given in Section 3.4 for separative algebras.

The following lemma follows immediately from the definitions of filter and ideal.

\begin{lm}\label{3.2.2} The intersection of any set of filters (ideals) is a filter (ideal).
\end{lm}

Let $A\subseteq X$. We define $\mu(A)$ - the {\em multiplicative closure of} $A$, by putting $\mu(A)$ to be the intersection of all filters containing $A$. By Lemma \ref{3.2.2}, $\mu(A)$ is the smallest filter containing $A$. Analogously, the intersection of all ideals containing $A$, denoted by $\alpha(A)$ and called the {\em additive closure of} $A$, is the smallest ideal containing $A$.

\begin{lm}\label{3.2.3} (i) \ $\mu(A)=\bigcup_{i=1}^{\infty}A^i$,

\smallskip

\noindent(i') \ $\alpha(A)=\bigcup_{i=1}^{\infty}iA$,

\smallskip

\noindent(ii) \ \ a) \ \ If $F$ is a filter then $F=\mu(F)$,

b)	\ \ If $A\subseteq B$ then $\mu(A)\subseteq\mu(B)$,

c)	\ \ $A\subseteq\mu(A)$,

d)	\ \ $\mu(\mu(A))=\mu(A)$,

e)	\ \ $\mu(A\cup B)=\mu(A)\cup\mu(A)\mu(B)\cup\mu(B)$; if $F$ and $G$ are filters then $\mu(F\cup G)=F\cup FG\cup G$; if $F$ is a filter and $a\in X$ then $\mu(F\cup a)=F\cup F.\mu(a)\cup\mu(a)$.

\smallskip

\noindent(ii') \ a) \ \ 	If $I$ is an ideal then $I=\alpha(I)$,

b)	\ \ If $A\subseteq B$ then $\alpha(A)\subseteq\alpha(B)$,

c) \ \ 	$A\subseteq\alpha(A)$,

d)	\ \ $\alpha(\alpha(A))=\alpha(A)$,

e)	\ \ $\alpha(A\cup B)=\alpha(A)\cup(\alpha(A)+\alpha(B))\cup\alpha(B)$;  if $I$ and $J$ are ideals
then $\alpha(I\cup J)=I\cup (I+J)\cup J$; if $I$ is an ideal and $a\in X$ then $\alpha(I\cup a)=I\cup I.\alpha(a)\cup\alpha(a)$.
\end{lm}

\doc (i) To prove the equality (i) it is enough to show that $\bigcup_{i=1}^{\infty}A^i$ is the smallest filter containing $A$. By Lemma \ref{3.1.2}(iv), we have  $(\bigcup_{i=1}^{\infty}A^i).(\bigcup_{i=1}^{\infty}A^i)\subseteq\bigcup_{i,j=1}^{\infty}A^i.A^j= \bigcup_{i,j=1}^{\infty}A^{i+j}\subseteq(\bigcup_{i=1}^{\infty}A^i)$,  so $\bigcup_{i=1}^{\infty}A^i$ is a filter, which obviously contains $A$. To prove that $\bigcup_{i=1}^{\infty}A^i$ is the smallest filter containing $A$, let $\alpha$ be a filter and $A\subseteq\alpha$. Applying Lemma \ref{3.1.2}(ii), we can show by induction on $i$ that $A^i\subseteq\alpha^i\subseteq\alpha$ and consequently $\bigcup_{i=1}^{\infty}A^i\subseteq\alpha$.

(i') can be shown similarly.

(ii) The proof of the conditions a),\ b),\ c) and d) follow directly from the definition of $\mu$. To prove condition e), we shall show that the set $F\cup FG\cup G$, where $F=\mu(A)$ and $G=\mu(B)$, is the smallest filter containing $A\cup B$.

By Lemma \ref{3.1.5}(ii), we obtain $$(F\cup FG\cup G)^2=F^2\cup F^2G\cup FG\cup F^2G^2\cup FG^2\cup G^2\subseteq F\cup FG\cup G.$$
This shows that $F\cup FG\cup G$ is a filter containing $F$ and $G$ and hence $A$ and $B$. To show that $F\cup FG\cup G$ is the smallest filter containing $A$ and $B$, let $\gamma$ be a filter such that $A\subseteq\gamma$ and $B\subseteq\gamma$, so we have $F\subseteq\gamma$ and $G\subseteq\gamma$. Then $F\cup G\subseteq\gamma$, $FG\subseteq\gamma\gamma\subseteq\gamma$ and consequently $F\cup FG\cup G\subseteq\gamma$.

The proof of (ii') can be obtained in a similar way. \sq

\begin{pro}\label{3.2.4}
Let $F$ be a filter and $I$ be an ideal. Then: \\
(i) $F-I$ is a filter, \\
(i') $\frac{I}{F}$ is an ideal, \\
(ii) If $I\cap(F-I)\not=\emptyset$ then $F\cap I\not=\emptyset$, \\
(iii) If $F\cap\frac{I}{F}\not=\emptyset$ then $F\cap I\not=\emptyset$, \\
(iv) If $(F-I)\cap\frac{I}{F}\not=\emptyset$ then $F\cap I\not=\emptyset$.
\end{pro}

\doc We shall proof (iv); the proofs of the other conditions are similar. Applying Proposition \ref{3.1.3}, we obtain:

$(F-I)\cap\frac{I}{F}\not=\emptyset \longleftrightarrow F\cap(I+\frac{I}{F})\not=\emptyset$; since  $I+\frac{I}{F}\subseteq\frac{I+I}{F}\subseteq\frac{I}{F}$, we get that $F\cap\frac{I}{F}\not=\emptyset$. \sq

\begin{lm}\label{3.2.5}
If $\mu(A)\cap\alpha(B)\not=\emptyset$ then there exist finite subsets $A'\subseteq A$ and $B'\subseteq B$ such that $\mu(A')\cap\alpha(B')\not=\emptyset$.
\end{lm}

\doc Let
\begin{equation}\label{(1)}
	\mu(A)\cap\alpha(B)\not=\emptyset.
\end{equation}
By Lemma \ref{3.2.3}(i)(i'), we have that
\begin{equation}\label{(2)}
	\mu(A)=\bigcup_{i=1}^{\infty}A^i\mbox{\ \ \ \ \ \    and}
\end{equation}
\begin{equation}\label{(3)}	
    \alpha(B)=\bigcup_{j=1}^{\infty}jB.
\end{equation}

From (\ref{(1)}), (\ref{(2)}) and (\ref{(3)}), we obtain that for some $x\in X$, $x\in\bigcup_{i=1}^{\infty}A^i$ and $x\in\bigcup_{j=1}^{\infty}jB$. Then for some $i$ and $j$ we have that
\begin{equation}\label{(4)}	
x\in A^i\mbox{\ \ \ \  \ \  and}
\end{equation}
\begin{equation}\label{(5)}
x\in jB.
\end{equation}

It follows from (\ref{(4)}) that there exist a set $A'=\{a_1,\dots,a_i\}\subseteq A$ such that $x\in \{a_1,\dots ,a_i\}$. From here we obtain that $\{a_1\dots a_i\}\subseteq\mu(A')$ and consequently
\begin{equation}\label{(6)}
x\in\mu(A')\subseteq(A).
\end{equation}

In an analogous way we obtain from (\ref{(5)}) that there exists a finite subset $B'=\{b_1,\dots,b_j\}\subseteq B$ such that
\begin{equation}\label{(7)}
x\in\alpha(B')\subseteq\alpha(B).
\end{equation}

Then from (\ref{(1)}) and (\ref{(6)}) and (\ref{(7)}) we obtain
\begin{equation}\label{(8)}	
\mu(A')\cap\alpha(B')\not=\emptyset
\end{equation}

Thus, for some finite subsets $A'\subseteq A$ and $B'\subseteq B$, we have $\mu(A')\cap\alpha(B')\not=\emptyset$. \sq

\subsection{Separative algebras}

Let $\underline{X}=(X, . , + )$ be a preseparative algebra. For $x,y\in X$ define $$x\leq y\ \mbox{ iff }\ \mu(x)\cap\alpha(y)\not=\emptyset.$$

\begin{defi}\label{3.3.1}
\rm A preseparative algebra $\underline{X}=(X, . , + )$ is called a {\em separative algebra} if the following axiom is satisfied:

\smallskip

\noindent$(Sep_{\, 0})$ The relation $\leq$ is transitive.

\smallskip

\noindent A separative algebra X is called a {\em convex space} if the operations $``."$ and $``+"$ coincide. In this case the filters and the ideals are called {\em convex sets} and the prime filters correspond to the notion of {\em half-space}.
\end{defi}

Convex  spaces  have  been  studied  by  several  authors: Tagamlitzki \cite{T}, Prodanov \cite{P4,P5}, Bair \cite{B}, Bryant \cite{Br}, Bryant and Webster \cite{BrW}.




We will now give several examples of separative algebras.

\begin{exa}\label{3.3.2}
\rm Let $\underline{L}=(L, \vee, \wedge, 0, 1)$ be a distributive lattice and for $x,y\in X$ define $x\times y=\{z\in L: z\geq x\wedge y\}$ and $x + y=\{z\in L: z\leq x\vee y\}$ (see Example \ref{ex2}). Then $\underline{X}$ is a separative algebra.
\end{exa}

\begin{exa}\label{3.3.3}
\rm Let $\underline{X}=(X, 1, +, . )$ be a commutative ring and for $x,y\in X$ define $x\times y=x.y$ and $x + y=A(x,y)$ , where $A(x,y)$ is the ring-ideal generated by the set $\{x,y\}$ (see Example \ref{ex1}). Then $\underline X$ is a separative algebra.
\end{exa}

\begin{exa}\label{3.3.4}
\rm Let $\underline{X}$ be a real linear space. For arbitrary $a,b\in X$, we set $a\times b=a + b=\{ta+(1-t) b: 0\leq t\leq 1\}$. Then $\underline{X}$ is a convex space.
\end{exa}

Apart from these starting examples, there is a number of other ones. It seems that whenever we have a satisfactory theory of prime ideals, then there is also a structure of separative algebra.

\begin{exa}\label{3.3.5}
\rm Let $X$ be an ordered linear topological space. Then $X$ is a separative algebra with respect to the operations

$a\times b=\{x\in X: \exists y\in ab\mbox{ with }x\leq y\}$,

$a + b=\{x\in X: \exists y\in ab\mbox{ with }x\geq y\}$,\\
where $ab=\{ta+(1-t) b: 0\leq t\leq 1\}$.
\end{exa}

\begin{exa}\label{3.3.6}
\rm Let $\underline{X}=(X, .)$ be a commutative semigroup. Then $\underline{X}$ is a convex space.
\end{exa}

The following lemma for filters and ideals is very important.

\begin{lm}\label{3.3.7} Let $\underline{X}$ be a separative algebra. Then for any $A,B\subseteq X$ and $x\in X$, we have that
if $\mu(A)\cap\alpha(B\cup x)\not=\emptyset$ and $\mu(x\cup A)\cap\alpha(B)\not=\emptyset$, then $\mu(A)\cap\alpha(B)\not=\emptyset$.
\end{lm}

\doc Suppose that the lemma does not hold and proceed to obtain a contradiction. Then for some $A,B\subseteq X$ and $x\in X$ we have that
\begin{equation}\label{1(1)}
\mu(A)\cap\alpha(B\cup x)\not=\emptyset,
\end{equation}

\begin{equation}\label{1(2)}
\mu(x\cup A)\cap\alpha(B)\not=\emptyset,\mbox{\ \  and}
\end{equation}

\begin{equation}\label{1(3)}
\mu(A)\cap\alpha(B)=\emptyset.
\end{equation}

By Lemma \ref{3.2.3}((ii)e)((ii')e), we obtain:

\begin{equation}\label{1(4)}
\mu(x\cup A)=\mu(A)\cup\mu(A)\mu(x)\cup\mu(x)\mbox{\ \  and}
\end{equation}

\begin{equation}\label{1(5)}
\alpha(B\cup x)=\alpha(B)\cup(\alpha(B)+\alpha(x))\cup\alpha(x).
\end{equation}

From (\ref{1(1)}), (\ref{1(3)}) and (\ref{1(5)}), we obtain that

\noindent either (a) $\mu(A)\cap(\alpha(B)+\alpha(x))\not=\emptyset$,

\noindent or \ \ \ \ \ (b) $\mu(A)\cap\alpha(x)\not=\emptyset$.

From (\ref{1(2)}), (\ref{1(3)}) and (\ref{1(4)}),  we obtain that

\noindent either (a') $(\mu(A)\mu(x))\cap\alpha(B)\not=\emptyset$,

\noindent or \ \ \ \ \ (b') $\mu(x)\cap\alpha(B)\not=\emptyset$.

So, we have to consider and to obtain a contradiction in each of the following combinations of cases:  (aa'), (ab'), (ba') and (bb'). As an example we shall treat of only the case (aa') - the remaining cases can be treated in a similar way. For the sake of brevity, we put $F=\mu(A)$, $I=\alpha(B)$; note that $F$ is a filter and $I$ is an ideal.
Now (a) and (a') become:

\noindent (a)	$F\cap(I+\alpha(x))\not=\emptyset$ and

\noindent (a') $I\cap(F.\mu(x))\not=\emptyset$.

Applying Proposition \ref{3.1.3} to (a) and (a'), we obtain

\begin{equation}\label{1(5')}
\mu(x)\cap\frac{I}{F}\not=\emptyset\mbox{ and}
\end{equation}

\begin{equation}\label{1(6)}	
\alpha(x)\cap(F-I)\not=\emptyset.
\end{equation}

By (\ref{1(5')}), we conclude that there exists $y\in X$ such that

\begin{equation}\label{1(7)}	
y\in\mu(x)\mbox{ and}
\end{equation}

\begin{equation}\label{1(8)}
y\in\frac{I}{F}.
\end{equation}

By (\ref{1(6)}), we obtain that for some $z\in X$ we have

\begin{equation}\label{1(9)}
z\in\alpha(x)\mbox{ and}
\end{equation}

\begin{equation}\label{1(10)}
z\in F-I.
\end{equation}

Conditions (\ref{1(7)}) and (\ref{1(9)}) are equivalent respectively to

\begin{equation}\label{1(11)}
y\cap\mu(x)\not=\emptyset\mbox{ and}
\end{equation}

\begin{equation}\label{1(12)}
z\cap\alpha(x)\not=\emptyset.
\end{equation}

Since $y\subseteq\alpha(y)$, using (\ref{1(11)}), we get

\begin{equation}\label{1(13)}
\mu(x)\cap\alpha(y)\not=\emptyset
\end{equation}
and, consequently, $x\leq y$.

Since $z\subseteq\mu(z)$, using (\ref{1(12)}), we get

\begin{equation}\label{1(14)}
\mu(z)\cap\alpha(x)\not=\emptyset
\end{equation}
and, consequently, $z\leq x$.

Now, by the axiom $(Sep_0)$, we obtain that $z\leq y$ and, consequently,

\begin{equation}\label{1(15)}
\mu(z)\cap\alpha(y)\not=\emptyset.
\end{equation}

By Proposition \ref{3.2.4}(i), $F-I$ is a filter and since, by (\ref{1(10)}), $z\in F-I$, we get that

\begin{equation}\label{1(16)}
\mu(z)\subseteq F-I.
\end{equation}

By Proposition \ref{3.2.4}(i'), $\frac{I}{F}$ is an ideal and since, by (\ref{1(8)}), $y\in\frac{I}{F}$, we get that

\begin{equation}\label{1(17)}
\alpha(y)\subseteq\frac{I}{F}.
\end{equation}

From (\ref{1(16)}) and (\ref{1(17)}), we get that

\begin{equation}\label{1(18)}
\mu(z)\cap\alpha(y)\subseteq(F-I)\cap\frac{I}{F}.
\end{equation}

By (\ref{1(15)}) and (\ref{1(18)}), we obtain that

\begin{equation}\label{1(19)}
(F-I)\cap\frac{I}{F}\not=\emptyset.
\end{equation}

Applying Proposition \ref{3.2.4}(iv), we obtain that $F\cap I\not=\emptyset$, i.e. $\mu(A)\cap\alpha(B)\not=\emptyset$, which contradicts (\ref{1(3)}). This completes the proof of the lemma. \sq

\begin{cor}\label{3.3.8}
If $F$ is a filter, $I$ is an ideal and $F\cap I=\emptyset$, then, for any $x\in X$, either $\mu(F\cup x)\cap I=\emptyset$ or $F\cap\alpha(I\cup x)=\emptyset$.
\end{cor}

\subsection{Separation theorem}

\begin{defi}\label{d3.4.1}
\rm Let $\underline{X}=(X, ., +)$ be a preseparative algebra. The
following statement is called the  {\em Separation principle  for}  $X$:

\medskip

\noindent(Sep) If $F_0$ is a filter, $I_0$  is an ideal and $F_0\cap I_0=\emptyset$ then there
exist a prime filter $F$ and a prime ideal $I$ such that $F_0\subseteq F$, $I_0\subseteq I$ and $F\cap I=\emptyset$.
\end{defi}

The main aim of this section is the following:

\begin{theorem}\label{3.4.1}
(Separation theorem for separative algebras) Let $\underline{X}=(X, . , + )$ be a separative algebra. Then $\underline{X}$ satisfies the Separation principle (Sep).
\end{theorem}

\doc Let $F_0$ be a filter in $\underline{X}$, $I_0$ be an ideal in $\underline{X}$ and $F_0\cap I_0=\emptyset$.

Let $M=\{F: F\mbox{ is a filter in }\underline{X}, F_0\subseteq F\mbox{ and }F\cap I_0=\emptyset\}$. It is easy to see that $M$ with the set-inclusion $\subseteq$ is an inductive set and hence, by the Zorn lemma, $M$ has a maximal element, say $F$.

Let $N=\{I: I\mbox{ is an ideal}, I_0\subseteq I\mbox{ and }F\cap I=\emptyset\}$. The set $N$ supplied with the set-inclusion is also an inductive set and hence, by the Zorn lemma, it has a maximal element, say $I$. We shall show that $F$ is a prime filter and  $I$ is a prime ideal.

Since $F$ is a filter, $I$ is an ideal and $F\cap I=\emptyset$, it is enough to show that $F\cup I=X$. Let $x\in X$. We shall show that either $x\in F$ or $x\in I$. Since $F\cap I=\emptyset$,  Corollary \ref{3.3.8} implies that either $\mu(F\cup x)\cap I=\emptyset$ or $F\cap\alpha(I\cup x)=\emptyset$.

\medskip

\noindent{\em Case 1:} $\mu(F\cup x)\cap I=\emptyset$. Since $I_0\subseteq I$, we obtain that $\mu(F\cup x)\cap I_0=\emptyset$. We  also have that $F_0\subseteq F\subseteq\mu(F\cup x)$. From here we obtain that the filter $\mu(F\cup x)\in M$. By the maximality of $F$ in $M$, we obtain that  $\mu(F\cup x)=F$, and hence  $x\in F$.

\medskip

\noindent{\em Case 2:} $F\cap\alpha(I\cup x)=\emptyset$. Since $I_0\subseteq I\subseteq\alpha(I\cup x)$, we obtain that $\alpha(I\cup x)\in N$. Then, by the maximality of $I$ in $N$, we obtain that $\alpha(I\cup x)=I$, and hence $x\in I$.

So we have found a prime filter $F\supseteq F_0$ and a prime ideal $I\supseteq I_0$ such that $F\cap I=\emptyset$, which proves the theorem. \sq


Let's note that Theorem \ref{3.4.1} generalizes a few well known statements: the Stone separation theorem for filters and ideals in distributive lattices \cite{Sto2} and in Boolean algebras \cite{Sto1}, as well as the separation theorem for convex sets in convex spaces from \cite{T}.

\begin{theorem}\label{3.4.2} Let $\underline{X}=(X, . , + )$ be a preseparative algebra. Then the following conditions are equivalent:

(i) $\underline{X}$ is a separative algebra,

(ii) $\underline{X}$ satisfies the Separation principle (Sep).
\end{theorem}

\doc The implication (i)$\longrightarrow$(ii) is just Theorem \ref{3.4.1}. For the converse implication (ii)$\longrightarrow$(i), we have to show that (Sep) implies $(Sep_0)$ (see Definition \ref{3.3.1} for $(Sep_0)$). So, let $a,b,c\in X$,

\begin{equation}\label{2(1)}	
a\leq b\ (\mbox{ i.e., }\mu(a)\cap\alpha(b)\not=\emptyset)\mbox{ and}
\end{equation}

\begin{equation}\label{2(2)}	
b\leq c\ (\mbox{ i.e., }\mu(b)\cap\alpha(c)\not=\emptyset)
\end{equation}

\noindent and suppose that

\begin{equation}\label{2(3)}	
a\nleqslant c\ (\mbox{ i.e., }\mu(a)\cap\alpha(c)=\emptyset).
\end{equation}
Then  (\ref{2(3)}) and (Sep) imply
that there exist a prime filter $F$ and and a prime ideal $I$ such
that

\begin{equation}\label{2(4)}	
F\cap I=\emptyset\ (\mbox{i.e. }X\setminus F=I),
\end{equation}

\begin{equation}\label{2(5)}	
\mu(a)\subseteq F\mbox{ and}
\end{equation}

\begin{equation}\label{2(6)}	
\alpha(c)\subseteq I.
\end{equation}

From (\ref{2(1)}) and (\ref{2(5)}) we obtain

\begin{equation}\label{2(7)}	
F\cap\alpha(b)\not=\emptyset.
\end{equation}

From (\ref{2(2)}) and (\ref{2(6)}) we obtain

\begin{equation}\label{2(8)}	
\mu(b)\cap I\not=\emptyset.
\end{equation}

For the element $b$ we have, by (\ref{2(4)}), that either $b\in F$ or $b\in I$.

\medskip

\noindent{\em Case 1:} $b\in F$. Then $\mu(b)\subseteq F$ and, by (\ref{2(8)}), we obtain that $F\cap I\not=\emptyset$ - a contradiction with (\ref{2(4)}).

\medskip

\noindent{\em Case 2:} $b\in I$. Then $\alpha(b)\subseteq I$ and, by (\ref{2(7)}), we obtain that $F\cap I\not=\emptyset$ - again a contradiction with (\ref{2(4)}).

This completes the proof of the theorem. \sq

We shall conclude this section by showing that the Separation theorem is equiva\-lent to the following statement, which is a generalization of the well known Wallman's lemma:

\begin{theorem}\label{3.4.3}
Let $\underline{X}=(X, ., + )$ be a preseparative algebra. Then the following conditions are equivalent:

\smallskip

\noindent(i) $\underline{X}$ is a separative algebra,

\smallskip

\noindent((ii) (Wallman's lemma) Let $M$ be a filter in $X$ and let, for any prime filter $F\supseteq M$, an element $x_{F}\in F$ be chosen. Then there exists a finite number of prime filters $F_i\supseteq M$, $i=1,\dots,n$, such that $M\cap\alpha(\{x_{F_1},\dots,x_{F_n}\})\not=\emptyset$.
\end{theorem}

\doc  (i)$\longrightarrow$(ii).\  Let $\underline{X}$ be a separative algebra and M be a
filter in $\underline{X}$. Denote by $N$ the set of all elements $x_{F}$, chosen as in the condition of the Wallman's lemma. Then $M\cap\alpha(N)\not=\emptyset$. To prove this suppose the contrary. Then there exists a prime filter $F\supseteq M$ such that $F\cap\alpha(N)=\emptyset$. But this is impossible because $x_{F}\in N\subseteq\alpha(N)$. So, $M\cap\alpha(N)\not=\emptyset$. Now,  by   Lemma  \ref{3.2.5}, there  exists   a   finite   subset
$\{x_{F_1},\dots,x_{F_n}\}\subseteq N$ such that $M\cap\alpha(\{x_{F_1},\dots,x_{F_n}\})\not=\emptyset$.

\medskip

\noindent(ii)$\longrightarrow$(i).\   Suppose  the  Wallman's  lemma.  We  shall  prove  the
Separation princi\-ple (Sep). Suppose, for the sake of contradiction, that (Sep) does not hold. Then, for some filter $F_0$ and some ideal $I_0$ such that $F_0\cap I_0=\emptyset$, we have that any prime filter $F$ extending $F_0$ has a non-empty intersection with $I_0$, i.e, there exists $x_{F}\in F\cap I_0$. Then,  by  the  Wallman  lemma,  there  exists  a  finite  set
$\{x_{F_1},\dots,x_{F_n}\}$  such that $F_0\cap\alpha(\{x_{F_1},\dots,x_{F_n}\})\not=\emptyset$.
But $\{x_{F_1},\dots,x_{F_n}\}\subseteq I_0$, so that $\alpha(\{x_{F_1},\dots,x_{F_n}\})\subseteq I_0$, which implies $F_0\cap I_0\not =\emptyset$ - a contradiction. \sq

\subsection{Standardization of the operations}

Here we shall consider two couples of natural operations in a given separative algebra.

Let $\underline{X}=(X, \otimes , \oplus)$ be a separative algebra and, for any $a,b\in X$, define the following two new multivalued operations, called {\em convex operations}:
$$a.b=\mu(\{a,b\})\mbox{ and }a+b=\alpha(\{a,b\})$$

\begin{theorem}\label{3.6.1}
If $\underline{X}$ is a separative algebra then it remains separative algebra with respect to its convex operations.
\end{theorem}

\doc The easy proof follows from the observation that the filters and ideals with respect to convex operations remain the same. \sq

Let $\underline{X}=(X, \otimes , \oplus )$ be a separative algebra. For any $A\subseteq X$, let $\mu_{\rho}(A)$ be the intersection of all prime filters containing $A$, and $\alpha_{\rho}(A)$ be the intersections of all prime ideals containing $A$. A subset $A$ of $X$ will be called a {\em radical filter} (resp., a {\em radical ideal}) if $\mu_{\rho}(A)=A$ (resp., $\alpha_{\rho}(A)=A$).

It follows from the Separation theorem that if $A$ is an ideal (resp. filter), then
$$\alpha_{\rho}(A)=\{x\in X:\ \mu(x)\cap A\not=\emptyset\},\mbox{ (resp., }\mu_{\rho}(A)=\{x\in X:\ \alpha(x)\cap A\not=\emptyset\}).$$

The following two new operations in X are called {\em radical operations}:
$$a.b=\mu_{\rho}(\{a,b\})\mbox{ and }a+b=\alpha_{\rho}(\{a,b\}),$$
where  $a,b\in X$.

\begin{theorem}\label{3.6.2}
If $\underline{X}=(X, \otimes , \oplus )$ is a separative algebra then it remains separative algebra with respect to its radical operations.
\end{theorem}

The proof follows from the observation that the filters and ideals with respect to the radical operations are the radical filters and radical ideals with respect to the initial operations, but the order $\leq$ do not change. To show this, note that $\mu_{\rho}(a)=\mu_{\rho}(\mu(a))$ and $\alpha_{\rho}(b)=\alpha_{\rho}(\alpha(b))$. Then, by the above observation, we have that
$$\mu_{\rho}(a)=\mu_{\rho}(\mu(a))=\{x\in X:\ \alpha(x)\cap\mu(a)\not=\emptyset\}=\{x\in X:\ a\leq x\}\mbox{ and}$$
$$\alpha_{\rho}(b)=\alpha_{\rho}(\alpha(b))=\{x\in X:\ \mu(x)\cap\alpha(b)\not=\emptyset\}=\{x\in X:\ x\leq b\}.$$

Then $\mu_{\rho}(a)\cap\alpha_{\rho}(b)\not=\emptyset$ iff $\exists x$: $a\leq x$ and $x\leq b$ iff $a\leq b$. \sq

\subsection{Canonical representation}

Let $\underline{X}$ be a separative algebra. Then $\underline{X}$ has a canonical representation $\varphi:X\longrightarrow L$ into a distributive lattice with the properties from Corollary \ref{th6}. Now $\varphi$ has some additional properties.

First of all, the inequality $a\leq b$ takes place if and only if $\varphi(a)\subseteq\varphi(b)$ . Therefore $\varphi(a)=\varphi(b)$ if and only if the radical ideals containing $a$ contain $b$ and conversely. If we do not distinguishing such points (which is natural, if we are interested only in radical ideals and filters), $\varphi$ becomes an embedding.

Now the operations from Corollary \ref{th6}(v) look in the following manner:
$$a.b=\{x\in X:\ \varphi(x)\leq\varphi(a)\vee\varphi(b)\}\mbox{ and }a+b=\{x\in X:\ \varphi(x)\geq\varphi(a)\wedge\varphi(b)\},$$ where $a.b$ and $a+b$ are the radical operations. In particular, if the initial operations coincide with radical ones, as it is in Example \ref{3.3.4}, we can get the separative structure of $X$ from suitable embedding of $X$ into a distributive lattice.

Now, let $\underline{X}$ be a ring with the separative structure from Example \ref{3.3.3}, and let $\varphi:X\longrightarrow L$ be the canonical representation. Then $L$ can be identified with the distributive lattice of all finitely generated radical ideals of $X$ (the whole $X$ is included), and, for arbitrary $a\in X$, the image $\varphi(a)$ is the radical ideal in $X$ generated by $a$.

\subsection{Topological version of the separation theorem}

\begin{defi}\label{3.8.1}
\rm We shall say that a preseparative algebra $\underline{X}=(X, ., + )$ is {\em topologi\-cal}, if $X$ is endowed with a topology such that the mappings $a.x$ and $a+x$ are lower semi-continuous, i.e., for every $a\in X$, the multi-valued maps $$\varphi_a:X\longrightarrow X, \ \ x\mapsto a+x,\ \ \mbox{ and }\ \ \psi_a:X\longrightarrow X, \ \ x\mapsto a.x,$$ are lower semi-continuous. Recall that a multi-valued map $f:X\longrightarrow Y$ between two topological spaces $X$ and $Y$ is said to be {\em lower semi-continuous}\/ if, for every open subset $U$ of $Y$, the set $f^{-1}(U)$ is open in $X$ (here, as usual, $$f^{-1}(U)=\{x\in X:\ f(x)\cap U\not=\emptyset\});$$ equivalently, $f$ is lower semi-continuous if, for every $x_0\in X$ and every open subset $U$ of $Y$ with $U\cap f(x_0)\not=\emptyset$, there exists a neighborhood $V$ of $x_0$  in $X$ such that $U\cap f(x)\not=\emptyset$, for every $x\in V$. For $a+x$, for example, this means that if $a,b\in X$ and $U$ is an open set with $(a+b)\cap U\not=\emptyset$, then there exists a neighborhood $V$ of $b$ such that $(a+x)\cap U\not=\emptyset$, for each $x\in V$.

A topological preseparative algebra will be called a {\em separative space} if, for each open filter $U$ in $X$, the conditions $\alpha(a)\cap U\not=\emptyset$ and $b\in\mu(a)$ imply $\alpha(b)\cap U\not=\emptyset$.

A separative space $\underline{X}=(X, ., + )$ is called a {\em topological convex space} if the operations $``."$ and  $``+"$ in $X$ coincide (see  \cite{P4}, \cite{P5}).
\end{defi}

Clearly, every separative algebra $X$ endowed with the discrete topology is a separati\-ve space, but there are also analytical examples. Now we shall only note that if $\underline{X}$ is a topological preseparative algebra such that the topology of $\underline{X}$ has a basis from open filters, then $X$ is a separative space.

The next statement, which we include here without proof, is a topological version of the Separation theorem.

\begin{theorem}\label{3.8.2}
Let $\underline{X}$ be a separative space, $I_0$ be an ideal in $X$ and $F_0$ be an open filter in $X$ such that $F_0\cap I_0=\emptyset$. Then there exist a closed prime ideal $I$ and an open prime filter $F$ in X such that
$F_0\subseteq F$, $I_0\subseteq I$ and $F\cap I=\emptyset$.
\end{theorem}

For a proof of Theorem \ref{3.8.2} for topological convex spaces see \cite{T}. We shall notice only one application of the theorem which uses the separative (not convex) structure: Example \ref{3.3.5} and Theorem \ref{3.8.2} give the classical separation theorem in ordered linear spaces, and, in particular, the general representation theorem of Kadison \cite{K}.

\baselineskip = 1.0\normalbaselineskip


 \end{document}